\documentclass{article}%
\usepackage{amssymb}
\usepackage{amsmath}
\usepackage{graphicx}
\usepackage{amsfonts}%
\providecommand{\U}[1]{\protect\rule{.1in}{.1in}}
\begin{document}

\title{Adjoint Functors and Heteromorphisms}
\author{David Ellerman}
\date{}
\maketitle

\begin{abstract}
Category theory has foundational importance because it provides conceptual
lenses to characterize what is important in mathematics. Originally the main
lenses were universal mapping properties and natural transformations. In
recent decades, the notion of adjoint functors has moved to center-stage as
category theory's primary tool to characterize what is important in
mathematics. Our focus here is to present a theory of adjoint functors. The
basis for the theory is laid by first showing that the object-to-object
\textquotedblleft heteromorphisms" between the objects of different categories
(e.g., insertion of generators as a set to group map) can be rigorously
treated within category theory. The heteromorphic theory shows that all
adjunctions arise from the birepresentations of \ the heteromorphisms\ between
the objects of different categories.

\end{abstract}
\tableofcontents

\section{The Importance of Adjoints}

Category theory is of foundational importance in mathematics but it is not
\textquotedblleft foundational\textquotedblright\ in the sense normally
claimed by set theory. It does not try to provide some basic objects (e.g.,
sets) from which other mathematical objects can be constructed. Instead, the
foundational role of category theory lies in providing conceptual lenses to
characterize what is universal and natural in mathematics.\footnote{For
summary statements, see \cite{aw:sml}, \cite{lanmarq:cinc}, or \cite{ell:cu}.}
Two of the most important concepts are universal mapping properties and
natural transformations. These two concepts are combined in the notion of
adjoint functors. In recent decades, adjoint functors have emerged as the
principal lens through which category theory plays out its foundational role
of characterizing what is important in mathematics.\footnote{Some familiarity
with basic category theory is assumed but the paper is written to be
accessible to a broader audience than specialists. Whenever possible, I follow
MacLane \cite{mac:cwm} on notation and terminology.}

The developers of category theory, Saunders MacLane and Samuel Eilenberg,
famously said that categories were defined in order to define functors, and
functors were defined in order to define \textit{natural} transformations.
Their original paper $\cite{eilmac:gte}$ was entitled not \textquotedblleft
General Theory of Categories\textquotedblright\ but \textit{General Theory of
Natural Equivalences}. Adjoints were defined more than a decade later by
Daniel Kan $\cite{kan:af}$ but the realization of their foundational
importance has steadily increased over time $\cite{law:adj, lam:her}$. Now it
would perhaps not be too much of an exaggeration to see categories, functors,
and natural transformations as the prelude to defining adjoint functors. As
Steven Awodey put it in his recent text:

\begin{quote}
{\footnotesize The notion of adjoint functor applies everything that we've
learned up to now to unify and subsume all the different universal mapping
properties that we have encountered, from free groups to limits to
exponentials. But more importantly, it also captures an important mathematical
phenomenon that is invisible without the lens of category theory. Indeed, I
will make the admittedly provocative claim that adjointness is a concept of
fundamental logical and mathematical importance that is not captured elsewhere
in mathematics. }$\cite{aw:cat}$
\end{quote}

\noindent Other category theorists have given similar testimonials.

\begin{quote}
{\footnotesize To some, including this writer, adjunction is the most
important concept in category theory. \cite[p. 6]{wood:ord}}

{\footnotesize The isolation and explication of the notion of
\textit{adjointness} is perhaps the most profound contribution that category
theory has made to the history of general mathematical ideas.\cite[p.
438]{Gold:topoi}}

{\footnotesize Nowadays, every user of category theory agrees that
[adjunction] is the concept which justifies the fundamental position of the
subject in mathematics. \cite[p. 367]{tay:pfm}}
\end{quote}

\noindent Given the importance of adjoint functors in category theory and in
mathematics as a whole, it would seem worthwhile to further investigate the
concept of an adjunction. Where do adjoint functors come from; how do they
arise? In this paper we will present a theory of adjoint functors to address
these questions \cite{ell:taf, ell:ae}. One might well ask: \textquotedblleft
Where could such a theory come from?\textquotedblright\ 

Category theory groups together in \textit{categories} the mathematical
objects with some common structure (e.g., sets, partially ordered sets,
groups, rings, and so forth) and the appropriate morphisms between such
objects. Since the morphisms are between objects of similar structure, they
are ordinarily called \textquotedblleft homomorphisms\textquotedblright\ or
just \textquotedblleft morphisms\textquotedblright\ for short. But there have
always been other morphisms which occur in mathematical practice that are
between objects with different structures (i.e., in different categories) such
as the insertion-of-generators map from a set to the free group on that set.
In order to contrast these morphisms with the homomorphisms between objects
within a category, they might be called \textit{heteromorphisms} or, more
colorfully, \textit{chimera morphisms} (since they have a tail in one category
and a head in another category). The usual machinery of category theory
(bifunctors, in particular) can be adapted to give a rigorous treatment of heteromorphisms.

With a precise notion of heteromorphisms in hand, it can then be seen that
adjoint functors arise as the functors giving the representations, using
homomorphisms \textit{within} each category, of the heteromorphisms
\textit{between} two categories. And, conversely, given a pair of adjoint
functors, then heteromorphisms can be defined between (isomorphic copies of)
the two categories so that the adjoints arise out of the representations of
those heteromorphisms. Hence this heteromorphic theory shows where adjoints
\textquotedblleft come from\textquotedblright\ or \textquotedblleft how they
arise.\textquotedblright\ It would seem that this theory showing the origin of
adjoint functors was not developed in the conventional treatment of category
theory since heteromorphisms, although present in mathematical practice, are
not part of the usual machinery of category theory.

\section{Overview of the Theory of Adjoints}

The cross-category object-to-object morphisms $c:x\Rightarrow a$, called
\textit{heteromorphisms} (\textit{hets} for short) or \textit{chimera
morphisms,} will be indicated by double arrows ($\Rightarrow$) rather than
single arrows ($\rightarrow$). The first question is how do heteromorphisms
compose with one another? But that is not necessary. Chimera do not need to
`mate' with other chimera to form a `species' or category; they only need to
mate with the intra-category morphisms on each side to form other chimera. The
appropriate mathematical machinery to describe that is the generalization of a
group acting on a set to a generalized monoid or category acting on a set
(where each element of the set has a \textquotedblleft
domain\textquotedblright\ and a \textquotedblleft codomain\textquotedblright%
\ to determine when composition is defined). In this case, it is two
categories acting on a set, one on the left and one on the right. Given a
chimera morphism $c:x\Rightarrow a$ from an object in a category $\mathbf{X}$
to an object in a category $\mathbf{A}$ and morphisms $h:x^{\prime}\rightarrow
x$ in $\mathbf{X}$ and $k:a\rightarrow a^{\prime}$ in $\mathbf{A}$, the
composition $ch:x^{\prime}\rightarrow x\Rightarrow a$ is another chimera
$x^{\prime}\Rightarrow a$ and the composition $kc:x\Rightarrow a\rightarrow
a^{\prime}$ is another chimera $x\Rightarrow a^{\prime}$ with the usual
identity, composition, and associativity properties. Such an action of two
categories acting on a set on the left and on the right is exactly described
by a bifunctor $\operatorname*{Het}:\mathbf{X}^{op}\times\mathbf{A}%
\rightarrow\mathbf{Set}$ where $\operatorname*{Het}(x,a)=\{x\Rightarrow a\}$
and where $\mathbf{Set}$ is the category of sets and set functions. Thus the
natural machinery to treat object-to-object chimera morphisms \textit{between}
categories are het-bifunctors $\operatorname*{Het}:\mathbf{X}^{op}%
\times\mathbf{A}\rightarrow\mathbf{Set}$ that generalize the hom-bifunctors
$\operatorname*{Hom}:\mathbf{X}^{op}\times\mathbf{X}\rightarrow\mathbf{Set}$
used to treat object-to-object morphisms \textit{within} a category.

How might the categorical properties of the heteromorphisms be expressed using
homomorphisms? Represent the het-bifunctors using hom-functors on the left and
on the right (see any category theory text such as \cite{mac:cwm} for
Alexander Grothendieck's notion of a \textit{representable functor}). Any
bifunctor $\mathcal{D}:\mathbf{X}^{op}\times\mathbf{A}\rightarrow\mathbf{Set}$
is \textit{represented on the left}\footnote{This terminology
\textquotedblleft represented on the left\textquotedblright\ or
\textquotedblleft on the right\textquotedblright\ is used to agree with the
terminology for left and right adjoints.} if for each $x$ in $\mathbf{X}$
there is an object $Fx$ in $\mathbf{A}$ and an isomorphism
$\operatorname*{Hom}_{\mathbf{A}}(Fx,a)\cong\mathcal{D}(x,a)$ natural in $a$.
It is a standard result that the assignment $x\mapsto Fx$ extends to a functor
$F$ and that the isomorphism is also natural in $x$. Similarly, $\mathcal{D}$
is \textit{represented on the right} if for each $a$ there is an object $Ga$
in $\mathbf{X}$ and an isomorphism $\mathcal{D}(x,a)\cong\operatorname*{Hom}%
_{\mathbf{X}}(x,Ga)$ natural in $x$. And similarly, the assignment $a\mapsto
Ga$ extends to a functor $G$ and that the isomorphism is also natural in $a$.

If a het-bifunctor $\operatorname*{Het}:\mathbf{X}^{op}\times\mathbf{A}%
\rightarrow\mathbf{Set}$ is represented on both the left and the right, then
we have two functors $F:\mathbf{X}\rightarrow\mathbf{A}$ and $G:\mathbf{A}%
\rightarrow\mathbf{X}$ and the isomorphisms are natural in $x$ and in $a$:

\begin{center}
$\operatorname*{Hom}_{\mathbf{A}}(Fx,a)\cong\operatorname*{Het}(x,a)\cong%
\operatorname*{Hom}_{\mathbf{X}}(x,Ga)$.
\end{center}

\noindent It only remains to drop out the middle term $\operatorname*{Het}%
(x,a)$ to arrive at the \textit{pas de deux} of the `official' definition of a
pair of adjoint functors which does not mention heteromorphisms.

While a birepresentation of a het-bifunctor gives rise to an adjunction, do
all adjunctions arise in this manner? To round out the theory, we give an
\textquotedblleft adjunction representation theorem\textquotedblright\ which
shows how, given any adjunction $F:\mathbf{X}\rightleftarrows\mathbf{A}$ $:G$,
heteromorphisms can be defined between (isomorphic copies of) the categories
$\mathbf{X}$ and $\mathbf{A}$ so that (isomorphic copies of) the adjoints
arise from the representations on the left and right of the het-bifunctor.
Given any set function $f:X\rightarrow A$ from the set $X$ to a set $A$, the
graph $\operatorname*{graph}\left(  f\right)  =\left\{  \left(  x,f(x)\right)
:x\in X\right\}  \subseteq X\times A$ of the function is set-isomorphic to the
domain of the function $X$. The embedding $x\longmapsto\left(  x,f\left(
x\right)  \right)  $ maps $X$ to the set-isomorphic copy of $X$, namely
$\operatorname*{graph}\left(  f\right)  \subseteq X\times A$. That isomorphism
generalizes to categories and to functors between categories. Given any
functor $F:\mathbf{X}\rightarrow\mathbf{A}$, the domain category $\mathbf{X}$
is embedded in the product category $\mathbf{X}\times\mathbf{A}$ by the
assignment $x\mapsto(x,Fx)$ to obtain the isomorphic copy $\widehat
{\mathbf{X}}$ (which can be considered as the graph of the functor $F$). Given
any other functor $G:\mathbf{A}\rightarrow\mathbf{X}$, the domain category
$\mathbf{A}$ is embedded in the product category by $a\mapsto(Ga,a)$ to yield
the isomorphic copy $\widehat{\mathbf{A}}$ (the graph of the functor $G$). If
the two functors are adjoints, then the properties of the adjunction can be
nicely expressed by the commutativity within the one category $\mathbf{X}%
\times\mathbf{A}$ of \textquotedblleft hom-pair adjunctive
squares\textquotedblright\ where morphisms are pairs of homomorphisms\ (in
contrast to a \textquotedblleft het adjunctive square\textquotedblright%
\ defined later).

\begin{center}
$%
\begin{array}
[c]{ccccc}
& (x,Fx) & \overset{(f,Ff)}{\longrightarrow} & (Ga,FGa) & \\
(\eta_{x},1_{Fx}) & \downarrow & \searrow^{\left(  f,g\right)  } & \downarrow
& (1_{Ga},\varepsilon_{a})\\
& (GFx,Fx) & \overset{(Gg,g)}{\longrightarrow} & (Ga,a) &
\end{array}
\nolinebreak$

Hom-pair adjunctive square
\end{center}

\noindent The main diagonal $(f,g)$ in a commutative hom-pair adjunctive
square pairs together maps that are images of one another in the adjunction
isomorphism $\operatorname*{Hom}_{\mathbf{A}}(Fx,a)\cong\operatorname*{Hom}%
_{\mathbf{X}}(x,Ga)$. If $f\in\operatorname*{Hom}_{\mathbf{X}}(x,Ga)$,
$g=f^{\ast}\in\operatorname*{Hom}_{\mathbf{A}}(Fx,a)$ is the corresponding
homomorphism on the other side of the isomorphism between hom-sets called its
\textit{adjoint transpose} (or later \textquotedblleft adjoint
correlate\textquotedblright) and similarly $f=g^{\ast}$. Since the maps on top
are always in $\widehat{\mathbf{X}}$ and the maps on the bottom are in
$\widehat{\mathbf{A}}$, the main diagonal pairs of maps (including the
vertical maps)---which are ordinary morphisms in the product category---have
all the categorical properties of heteromorphisms from objects in
$\mathbf{X}\cong\widehat{\mathbf{X}}$ to objects in $\mathbf{A}\cong%
\widehat{\mathbf{A}}$. Hence the heteromorphisms are abstractly defined as the
\textit{pairs of adjoint transposes}, $\operatorname*{Het}%
(x,a)=\{(x,Fx)\overset{(f,f^{\ast})}{\longrightarrow}(Ga,a)\}$, and the
adjunction representation theorem is that (isomorphic copies of) the original
adjoints $F$ and $G$ arise from the representations on the left and right of
this het-bifunctor.

Heteromorphisms are formally treated using bifunctors of the form
$\operatorname*{Het}:\mathbf{X}^{op}\times\mathbf{A}\rightarrow\mathbf{Set}$.
Such bifunctors and generalizations replacing $\mathbf{Set}$ by other
categories have been studied by the Australian school under the name of
\textit{profunctors }\cite{kelly:enrich}, by the French school under the name
of \textit{distributors} \cite{ben:distr}, and by William Lawvere under the
name of \textit{bimodules} \cite{law:met}. However, the guiding interpretation
has been interestingly different. \textquotedblleft Roughly speaking, a
distributor is to a functor what a relation is to a mapping\textquotedblright%
\ \cite[p. 308]{bor:hca1} (and hence the name \textquotedblleft
profunctor\textquotedblright\ in the Australian school). For instance, if
$\mathbf{Set}$ was replaced by $\mathbf{2}$, then the bifunctor would just be
the characteristic function of a relation from $\mathbf{X}$ to $\mathbf{A}$.
Hence in the general context of enriched category theory, a \textquotedblleft
bimodule\textquotedblright\ $Y^{op}\otimes X\overset{\varphi}{\longrightarrow
\mathcal{V}}$ would be interpreted as a \textquotedblleft$\mathcal{V}$-valued
relation\textquotedblright\ and an element of $\varphi(y,x)$ would be
interpreted as the \textquotedblleft truth-value of the $\varphi$-relatedness
of $y$ to $x$\textquotedblright\ \cite[p. 158 (p. 28 of reprint)]{law:met}.

The subsequent development of profunctors-distributors-bimodules has been
along the lines suggested by that guiding interpretation. For instance,
composition is defined between distributors as \textquotedblleft
relational\textquotedblright\ generalizations of functors to define a category
of distributors in analogy with composition defined between relations as
generalizations of functions which allows the definition of a category of
relations \cite[Chapter 7]{bor:hca1}.

The heteromorphic interpretation of the bifunctors $\mathbf{X}^{op}%
\times\mathbf{A}\rightarrow\mathbf{Set}$ is rather different. Each such
bifunctor is taken as defining how the chimeras in $\operatorname*{Het}(x,a)$
compose with morphisms in $\mathbf{A}$ on one side and with morphisms in
$\mathbf{X}$ on the other side to form other chimeras. This provides the
formal treatment of the heteromorphisms that have always existed in
mathematical practice. The principal novelty here is the use of the chimera
morphism interpretation of these bifunctors to carry out a whole program of
interpretation for adjunctions, i.e., a \textit{theory}\ of adjoint functors.
In the concrete examples, heteromorphisms have to be \textquotedblleft
found\textquotedblright\ as is done in the broad classes of examples treated
here. However, in general, the adjunction representation theorem shows how
abstract heteromorphisms (pairs of adjoint transposes in the product category
$\mathbf{X}\times\mathbf{A}$) can always be found so that any adjunction
arises (up to isomorphism) out of the representations on the left and right of
the het-bifunctor of such heteromorphisms. Following this summary, we now turn
to a slower development of the theory along with examples.

\section{The Heteromorphic Theory of Adjoints}

\subsection{Definition and Directionality of Adjoints}

There are many equivalent definitions of adjoint functors $\cite{mac:cwm}$,
but the most `official' one seems to be the one using a natural isomorphism of
hom-sets. Let $\mathbf{X}$ and $\mathbf{A}$ be categories and
$F:\mathbf{X\rightarrow A}$ and $G:\mathbf{A\rightarrow X}$ functors between
them. Then $F$ and $G$ are said to be a pair of \textit{adjoint functors} or
an \textit{adjunction}, written $F\dashv G$, if for any $x$ in $\mathbf{X}$
and $a$ in $\mathbf{A}$, there is an isomorphism $\phi$ natural in $x$ and in
$a$:

\begin{center}
${\phi}_{x,a}:\operatorname*{Hom}_{\mathbf{A}}(Fx,a)\;{\cong\;}%
\operatorname*{Hom}_{\mathbf{X}}(x,Ga)$.
\end{center}

\noindent With this standard way of writing the isomorphism of hom-sets, the
functor $F$ on the left is called the \textit{left adjoint} and the functor
$G$ on the right is the \textit{right adjoint}. Maps associated with each
other by the adjunction isomorphism ("adjoint transposes" of one another) are
indicated by an asterisk so if $g:Fx\rightarrow a$ then $g^{\ast}:x\rightarrow
Ga$ is the associated map ${\phi}_{x,a}(g)=g^{\ast}$ and similarly if
$f:x\rightarrow Ga$ then $\phi_{x,a}^{-1}(f)=f^{\ast}:Fx\rightarrow a$ is the
associated map.

In much of the literature, adjoints are presented in a seemingly symmetrical
fashion so that there appears to be no directionality of the adjoints between
the categories $\mathbf{X}$ and $\mathbf{A}$. But there is a directionality
and it is important in understanding adjoints. Both the maps that appear in
the adjunction isomorphism, $Fx\rightarrow a$ and $x\rightarrow Ga$, go from
the \textquotedblleft$x$-thing\textquotedblright\ (i.e., either $x$ or the
image $Fx$) to the \textquotedblleft$a$-thing\textquotedblright\ (either the
image $Ga$ or $a$ itself), so we see a direction emerging from $\mathbf{X}$ to
$\mathbf{A}$. That direction of an adjunction is the direction of the left
adjoint (which goes from $\mathbf{X}$ to $\mathbf{A}$). Then $\mathbf{X}$
might called the \textit{sending} category and $\mathbf{A}$ the
\textit{receiving} category.\footnote{Sometimes adjunctions are written with
this direction as in the notation $\langle F,G,{\phi\rangle}:\mathbf{X}%
{\rightharpoonup}$ $\mathbf{A}$ (MacLane $\cite[p. 78]{mac:cwm}$). This also
allows the composition of adjoints to be defined in a straightforward manner
(MacLane $\cite[p. 101]{mac:cwm}$).}

In the theory of adjoints presented here, the directionality of adjoints
results from being representations of heteromorphisms\ which have that
directionality. \ Such morphisms can exhibited in concrete examples of
adjoints (see the later examples). To abstractly define chimera morphisms\ or
heteromorphisms that work for all adjunctions, we turn to the presentation of
adjoints using adjunctive squares.

\subsection{Embedding Adjunctions in a Product Category}

Our approach to a theory of adjoints uses a certain \textquotedblleft
adjunctive square\textquotedblright\ diagram that is in the product category
$\mathbf{X}\times\mathbf{A}$ associated with an adjunction $F:\mathbf{X}%
\rightleftarrows\mathbf{A}:G$. With each object $x$ in the category
$\mathbf{X}$, we associate the element $\widehat{x}=(x,Fx)$ in the product
category $\mathbf{X}\times\mathbf{A}$ so that $Ga$ would have associated with
it $\widehat{Ga}=(Ga,FGa)$.\ With each morphism in $\mathbf{X}$ with the form
$h:x^{\prime}\rightarrow x$, we associate the morphism $\widehat
{h}=(h,Fh):\widehat{x^{\prime}}=(x^{\prime},Fx^{\prime})\rightarrow\widehat
{x}=(x,Fx)$ in the product category $\mathbf{X}\times\mathbf{A}$ (maps compose
and diagrams commute component-wise). Thus the mapping of $x$ to $(x,Fx)$
extends to an embedding $(1_{\mathbf{X}},F):\mathbf{X}\rightarrow
\mathbf{X}\times\mathbf{A}$ whose image $\widehat{\mathbf{X}}$ (the graph of
$F$) is isomorphic with $\mathbf{X}$.

With each object $a$ in the category $\mathbf{A}$, we associate the element
$\widehat{a}=(Ga,a)$ in the product category $\mathbf{X}\times\mathbf{A}$ so
that $Fx$ would have associated with it $\widehat{Fx}=(GFx,Fx)$. With each
morphism in $\mathbf{A}$ with the form $k:a\rightarrow a^{\prime},$we
associate the morphism $\widehat{k}=(Gk,k):(Ga,a)\rightarrow(Ga^{\prime
},a^{\prime})$ in the product category $\mathbf{X}\times\mathbf{A}$. The
mapping of $a$ to $(Ga,a)$ extends to an embedding $(G,1_{\mathbf{A}%
}):\mathbf{A}\rightarrow\mathbf{X}\times\mathbf{A}$ whose image $\widehat
{\mathbf{A}}$ (the graph of $G$) is isomorphic to $\mathbf{A}$.

The adjoint transpose of the identity map $1_{Fx}\in\operatorname*{Hom}%
_{A}\left(  Fx,Fx\right)  $ is the unit morphism $\eta_{x}:x\rightarrow
GFx\in\operatorname*{Hom}_{X}(x,GFx)$. That pair $(\eta_{x},1_{Fx}%
):(x,Fx)\rightarrow(GFx,Fx)$ of adjoint transposes is the left vertical
`heteromorphism' in the hom-pairs adjunctive square diagram. We use the
raised-eyebrow quotes on `heteromorphism' since it is a perfectly ordinary
homomorphism in the product category $\mathbf{X}\times\mathbf{A}$ which plays
the role of a heteromorphism from $\widehat{\mathbf{X}}$, the isomorphic copy
of $\mathbf{X}$, to $\widehat{\mathbf{A}}$, the isomorphic copy of
$\mathbf{A}$, both subcategories of $\mathbf{X}\times\mathbf{A}$. The adjoint
transpose of the identity map $1_{Ga}\in\operatorname*{Hom}_{X}(Ga,Ga)$ is the
counit morphism $\varepsilon_{a}:FGa\rightarrow a\in\operatorname*{Hom}%
_{A}(FGa,a)$. That pair $(1_{G},\varepsilon_{a}):(Ga,FGa)\rightarrow(Ga,a)$ of
adjoint transposes is the right vertical `heteromorphism' in the adjunctive
square diagram.

These various parts can then be collected together in the (hom-pair adjunctive
square diagram of the representation theorem.

\begin{center}
$%
\begin{array}
[c]{rcccl}
& (x,Fx) & \overset{(f,Ff)}{\rightarrow} & (Ga,FGa) & \\
(\eta_{x},1_{Fx}) & \downarrow & \searrow^{\left(  f,g\right)  } & \downarrow
& (1_{Ga},\varepsilon_{a})\\
& (GFx,Fx) & \overset{(Gg,g)}{\longrightarrow} & (Ga,a) &
\end{array}
\medskip\nolinebreak$

Hom-pair Adjunctive Square Diagram
\end{center}

The adjunctive square diagram conveniently represents the properties of an
adjunction in the format of commutative squares. \ The map on the top is in
$\widehat{\mathbf{X}}$ and the map on the bottom is in $\widehat{\mathbf{A}}$
and the vertical maps as well as the main diagonal $(f,g)$ in a commutative
adjunctive square are morphisms from $\widehat{\mathbf{X}}$-objects to
$\widehat{\mathbf{A}}$-objects.

Given $f:x\rightarrow Ga$, the rest of the diagram is determined by the
requirement that the square commutes. Commutativity in the second component
uniquely determines that $g=g1_{Fx}=\varepsilon_{a}Ff$ so $g=f^{\ast
}=\varepsilon_{a}Ff$ \ is the map associated with $f$ in the adjunction
isomorphism. \ Commutativity in the first component is the universal mapping
property (UMP) factorization of any given $f:x\rightarrow Ga$ through the unit
$x\overset{\eta_{x}}{\longrightarrow}GFx\overset{Gf^{\ast}}{\longrightarrow
}Ga=x\overset{f}{\longrightarrow}Ga$ which is often pictured as:

\begin{center}
$%
\begin{array}
[c]{cccc}
& x &  & \\
\eta_{x} & \downarrow & \searrow^{f} & \\
& GFx & \overset{Gf^{\ast}}{\longrightarrow} & Ga\\
& Fx & \overset{\exists!f^{\ast}}{\longrightarrow} & a
\end{array}
\nolinebreak$

Hom-pair adjunctive square south-west of the diagonal.
\end{center}

Similarly, if we were given $g:Fx\rightarrow a,$ then commutativity in the
first component implies that $f=1_{Ga}f=Gg\eta_{x}=g^{\ast}$. And
commutativity in the second component is the UMP factorization of any given
$g:Fx\rightarrow a$ through the counit $Fx\overset{Fg^{\ast}}{\longrightarrow
}FGa\overset{\varepsilon_{a}}{\longrightarrow}a=Fx\overset{g}{\longrightarrow
}a$ which is usually pictured as:

\begin{center}
$%
\begin{array}
[c]{cccc}%
x & \overset{\exists!g^{\ast}}{\longrightarrow} & Ga & \\
Fx & \overset{Fg^{\ast}}{\longrightarrow} & FGa & \\
& \searrow^{g} & \downarrow & \varepsilon_{a}\\
&  & a &
\end{array}
\nolinebreak$

Hom-pair adjunctive square north-east of the diagonal.
\end{center}

\noindent Splicing together the two triangles along the diagonals so that the
two diagonals form the hom-pair $\left(  f,g\right)  $ (and supplying the
identity maps $1_{Fx}$ and $1_{Ga}$ as required to form the left and right
vertical hom-pairs), the hom-pair adjunctive square is put back together.

\subsection{Heteromorphisms and Het-bifunctors}

Heteromorphisms\ (in contrast to homomorphisms) are like mongrels or chimeras
that do not fit into either of the two categories. Since the inter-category
heteromorphisms are not morphisms in either of the categories, what can we say
about them? The one thing we can reasonably say is that heteromorphisms can be
precomposed or postcomposed with morphisms within the categories (i.e.,
intra-category morphisms) to obtain other heteromorphisms.\footnote{The
chimera genes are dominant in these mongrel matings. While mules cannot mate
with mules, it is `as if' mules could mate with either horses or donkeys to
produce other mules.} This is easily formalized using bifunctors similar to
the hom-bifunctors $\operatorname*{Hom}(x,y)$ in homomorphisms within a
category. Using the sets-to-groups example to guide intuition, one might think
of $\operatorname*{Het}(x,a)=\{x\overset{c}{\Rightarrow}a\}$ as the set of set
functions from a set $x$ to a group $a$. For any $\mathbf{A}$-morphism
$k:a\rightarrow a^{\prime}$ and any chimera morphism $x\overset{c}%
{\Rightarrow}a$, intuitively there is a composite chimera morphism
$x\overset{c}{\Rightarrow}a\overset{k}{\rightarrow}a^{\prime}=x\overset
{kc}{\Rightarrow}a^{\prime}$, i.e., $k$ induces a map $\operatorname*{Het}%
(x,k):\operatorname*{Het}(x,a)\rightarrow\operatorname*{Het}(x,a^{\prime})$.
For any $\mathbf{X}$-morphism $h:x^{\prime}\rightarrow x$ and chimera morphism
$x\overset{c}{\Rightarrow}a$, intuitively there is the composite chimera
morphism $x^{\prime}\overset{h}{\rightarrow}x\overset{c}{\Rightarrow
}a=x^{\prime}\overset{ch}{\Rightarrow}a$, i.e., $h$ induces a map
$\operatorname*{Het}(h,a):\operatorname*{Het}(x,a)\rightarrow
\operatorname*{Het}(x^{\prime},a)$ (note the reversal of direction). The
induced maps would respect identity and composite morphisms in each category.
Moreover, composition is associative in the sense that $(kc)h=k(ch)$. This
means that the assignments of sets of chimera morphisms $\operatorname*{Het}%
(x,a)=\{x\overset{c}{\Rightarrow}a\}$ and the induced maps between them
constitute a \textit{bifunctor} $\operatorname*{Het}:\mathbf{X}^{op}%
\times\mathbf{A}\rightarrow\mathbf{Set}$ (contravariant in the first variable
and covariant in the second).

With this motivation, we may turn around and define \textit{heteromorphisms}%
\ from $\mathbf{X}$-objects to $\mathbf{A}$-objects as the elements in the
values of a bifunctor $\operatorname*{Het}:\mathbf{X}^{op}\times
\mathbf{A}\rightarrow\mathbf{Set}$. This would be analogous to defining the
homomorphisms in $\mathbf{X}$ as the elements in the values of a given
hom-bifunctor $\operatorname*{Hom}_{\mathbf{X}}:\mathbf{X}^{op}\times
\mathbf{X}\rightarrow\mathbf{Set}$ and similarly for $\operatorname*{Hom}%
_{\mathbf{A}}:\mathbf{A}^{op}\times\mathbf{A}\rightarrow\mathbf{Set}$.

With heteromorphisms rigorously described using het-bifunctors, we can use
Grothendieck's notion of a representable functor to show that an adjunction
arises from a het-bifunctor $\operatorname*{Het}:\mathbf{X}^{op}%
\times\mathbf{A}\rightarrow\mathbf{Set}$ that is \textquotedblleft
birepresentable\textquotedblright\ in the sense of being representable on both
the left and right.

Given any bifunctor $\operatorname*{Het}:\mathbf{X}^{op}\times\mathbf{A}%
\rightarrow\mathbf{Set}$, it is \textit{representable on the left} if for each
$\mathbf{X}$-object $x$, there is an $\mathbf{A}$-object $Fx$ that represents
the functor $\operatorname*{Het}(x,-)$, i.e., there is an isomorphism
$\psi_{x,a}:\operatorname*{Hom}_{\mathbf{A}}(Fx,a)\cong\operatorname*{Het}%
(x,a)$ natural in $a$. For each $x$, let $h_{x}$ be the image of the identity
on $Fx$, i.e., $\psi_{x,Fx}(1_{Fx})=h_{x}\in\operatorname*{Het}(x,Fx)$. We
first show that $h_{x}$ is a universal element for the functor
$\operatorname*{Het}(x,-)$ and then use that to complete the construction of
$F$ as a functor. For any $c\in\operatorname*{Het}(x,a)$, let $g(c)=\psi
_{x,a}^{-1}(c):Fx\rightarrow a$. Then naturality in $a$ means that the
following diagram commutes.

\begin{center}
$%
\begin{array}
[c]{ccccc}
& \operatorname*{Hom}_{\mathbf{A}}(Fx,Fx) & \cong & \operatorname*{Het}%
(x,Fx) & \\
_{\operatorname*{Hom}(Fx,g(c))} & \downarrow &  & \downarrow &
_{\operatorname*{Het}(x,g(c))}\\
& \operatorname*{Hom}_{\mathbf{A}}(Fx,a) & \cong & \operatorname*{Het}(x,a) &
\end{array}
\nolinebreak\medskip$

$\operatorname*{Het}\left(  x,a\right)  $ representable on the left
\end{center}

\noindent Chasing $1_{Fx}$ around the diagram yields that
$c=\operatorname*{Het}(x,g(c))(h_{x})$ which can be written as $c=g(c)h_{x}$.
Since the horizontal maps are isomorphisms, $g(c)$ is the unique map
$g:Fx\rightarrow a$ such that $c=gh_{x}$. Then $(Fx,h_{x})$ is a
\textit{universal element} (in MacLane's sense \cite[p. 57]{mac:cwm}) for the
functor $\operatorname*{Het}(x,-)$ or equivalently $1\overset{h_{x}%
}{\longrightarrow}\operatorname*{Het}(x,Fx)$ is a \textit{universal arrow}
\cite[p. 58]{mac:cwm} from $1$ (the one point set) to $\operatorname*{Het}%
(x,-)$. Then for any $\mathbf{X}$-morphism $j:x\rightarrow x^{\prime}$,
$Fj:Fx\rightarrow Fx^{\prime}$ is the unique $\mathbf{A}$-morphism such that
$\operatorname*{Het}(x,Fj)$ fills in the right vertical arrow in the following diagram.

\begin{center}
$%
\begin{array}
[c]{rcccc}
& 1 & \overset{h_{x}}{\longrightarrow} & \operatorname*{Het}(x,Fx) & \\
_{h_{x^{\prime}}} & \downarrow &  & \downarrow & _{\operatorname*{Het}%
(x,Fj)}\\
& \operatorname*{Het}(x^{\prime},Fx^{\prime}) & \overset{\operatorname*{Het}%
(j,Fx^{\prime})}{\longrightarrow} & \operatorname*{Het}(x,Fx^{\prime}) &
\end{array}
$
\end{center}

\noindent It is easily checked that such a definition of $Fj:Fx\rightarrow
Fx^{\prime}$ preserves identities and composition using the functoriality of
$\operatorname*{Het}(x,-)$ so we have a functor $F:\mathbf{X}\rightarrow
\mathbf{A}$. It is a further standard result that the isomorphism is also
natural in $x$ (e.g., \cite[p. 81]{mac:cwm} or the "parameter theorem"
\cite[p. 525]{macbirk:alg}).

Given a bifunctor $\operatorname*{Het}:\mathbf{X}^{op}\times\mathbf{A}%
\rightarrow\mathbf{Set}$, it is \textit{representable on the right} if for
each $\mathbf{A}$-object $a$, there is an $\mathbf{X}$-object $Ga$ that
represents the functor $\operatorname*{Het}(-,a)$, i.e., there is an
isomorphism $\varphi_{x,a}:\operatorname*{Het}(x,a)\cong\operatorname*{Hom}%
_{X}(x,Ga)$ natural in $x$. For each $a$, let $e_{a}$ be the inverse image of
the identity on $Ga$, i.e., $\varphi_{Ga,a}^{-1}(1_{Ga})=e_{a}\in
\operatorname*{Het}(Ga,a)$. For any $c\in\operatorname*{Het}(x,a)$, let
$f(c)=\varphi_{x,a}(c):x\rightarrow Ga$. Then naturality in $x$ means that the
following diagram commutes.

\begin{center}
$%
\begin{array}
[c]{ccccc}
& \operatorname*{Het}(Ga,a) & \cong & \operatorname*{Hom}_{\mathbf{X}%
}(Ga,Ga) & \\
_{\operatorname*{Het}(f(c),a)} & \downarrow &  & \downarrow &
_{\operatorname*{Hom}(f(c),Ga)}\\
& \operatorname*{Het}(x,a) & \cong & \operatorname*{Hom}_{\mathbf{X}}(x,Ga) &
\end{array}
\nolinebreak\medskip$

$\operatorname*{Het}\left(  x,a\right)  $ representable on the right
\end{center}

\noindent Chasing $1_{Ga}$ around the diagram yields that
$c=\operatorname*{Het}(f(c),a)(e_{a})=e_{a}f(c)$ so $(Ga,e_{a})$ is a
universal element for the functor $\operatorname*{Het}(-,a)$ and that
$1\overset{e_{a}}{\longrightarrow}\operatorname*{Het}(Ga,a)$ is a universal
arrow from $1$ to $\operatorname*{Het}(-,a)$. Then for any $\mathbf{A}%
$-morphism $k:a^{\prime}\rightarrow a$, $Gk:Ga^{\prime}\rightarrow Ga$ is the
unique $\mathbf{X}$-morphism such that $\operatorname*{Het}(Gk,a)$ fills in
the right vertical arrow in the following diagram.

\begin{center}
$%
\begin{array}
[c]{ccccc}
& 1 & \overset{e_{a}}{\longrightarrow} & \operatorname*{Het}(Ga,a) & \\
_{e_{a^{\prime}}} & \downarrow &  & \downarrow & _{\operatorname*{Het}%
(Gk,a)}\\
& \operatorname*{Het}(Ga^{\prime},a^{\prime}) & \overset{\operatorname*{Het}%
(Ga^{\prime},k)}{\longrightarrow} & \operatorname*{Het}(Ga^{\prime},a) &
\end{array}
$
\end{center}

\noindent In a similar manner, it is easily checked that the functoriality of
$G$ follows from the functoriality of $\operatorname*{Het}(-,a)$. Thus we have
a functor $G:\mathbf{A}\rightarrow\mathbf{X}$ such that $Ga$ represents the
functor $\operatorname*{Het}(-,a)$, i.e., there is a natural isomorphism
$\varphi_{x,a}:\operatorname*{Het}(x,a)\cong\operatorname*{Hom}_{\mathbf{X}%
}(x,Ga)$ natural in $x$. And in a similar manner, it can be shown that the
isomorphism is natural in both variables.

Thus given a bifunctor $\operatorname*{Het}:\mathbf{X}^{op}\times
\mathbf{A}\rightarrow\mathbf{Set}$ representable on \textit{both} sides, we
have the adjunction natural isomorphisms:

\begin{center}
$\operatorname*{Hom}_{\mathbf{A}}(Fx,a)\cong\operatorname*{Het}(x,a)\cong%
\operatorname*{Hom}_{\mathbf{X}}(x,Ga)$.
\end{center}

\noindent Starting with $c\in\operatorname*{Het}(x,a)$, the corresponding
$f\left(  c\right)  \in\operatorname*{Hom}_{\mathbf{X}}(x,Ga)$ and $g\left(
c\right)  \in\operatorname*{Hom}_{\mathbf{A}}(Fx,a)$ are called
\textit{adjoint correlates} of one another. Starting with $1_{Fx}%
\in\operatorname*{Hom}_{\mathbf{A}}(Fx,Fx)$, its adjoint correlates are the
\textit{het unit} $h_{x}\in\operatorname*{Het}(x,Fx)$ and the ordinary unit
$\eta_{x}\in\operatorname*{Hom}_{\mathbf{X}}(x,GFx)$ where this usual unit
$\eta_{x}$ might also be called the \textquotedblleft hom
unit\textquotedblright\ to distinguish it from its het correlate. Starting
with $1_{Ga}\in\operatorname*{Hom}_{\mathbf{X}}(Ga,Ga)$, its adjoint
correlates are the \textit{het counit} $e_{a}\in\operatorname*{Het}(Ga,a)$ and
the usual (hom) counit $\varepsilon_{a}\in\operatorname*{Hom}_{\mathbf{A}%
}(FGa,a)$. Starting with any $c\in\operatorname*{Het}\left(  x,a\right)  $,
the two factorizations $g(c)h_{x}=c=e_{a}f(c)$ combine to give what we will
later call the \textquotedblleft het adjunctive square\textquotedblright\ with
$c$ as the main diagonal [as opposed to the hom-pair adjunctive square
previously constructed which had $(f(c),g\left(  c\right)  )$ as the main diagonal].

There are cases (see below) where the het-bifunctor is only representable on
the left $\operatorname*{Hom}_{\mathbf{A}}(Fx,a)\cong\operatorname*{Het}(x,a)$
or on the right $\operatorname*{Het}(x,a)\cong\operatorname*{Hom}_{\mathbf{X}%
}(x,Ga)$, and in that case, it would make perfectly good sense to respectively
take $F:\mathbf{X}\rightarrow\mathbf{A}$ as a \textit{left half-adjunction} or
$G:\mathbf{A}\rightarrow\mathbf{X}$ as a \textit{right half-adjunction}. A
half-adjunction is the simplest expression of a universal mapping property,
and, of course, a left half-adjunction plus a right half-adjunction equals an
adjunction.\footnote{When there is a full adjunction, often only one
half-adjunction is important while the other half-adjunction is a rather
trivial piece of conceptual bookkeeping to round out the whole adjunction. In
the later examples, for the free-group/underlying-set adjunction, the left
half-adjunction carries the weight while in the Cartesian product adjunction,
it is the right half-adjunction.}

\subsection{Adjunction Representation Theorem}

Adjunctions may be and usually are presented without any thought to any
underlying heteromorphisms. \ However, given any adjunction, there is always
an \textquotedblleft abstract\textquotedblright\ associated \ het-bifunctor
given by the main diagonal maps in the commutative hom-pair adjunctive squares:

\begin{center}
$\operatorname*{Het}(\widehat{x},\widehat{a})=\{\widehat{x}=(x,Fx)\overset
{(f,f^{\ast})}{\longrightarrow}(Ga,a)=\widehat{a}\}$

Het-bifunctor for any adjunction from hom-pair adjunctive squares.
\end{center}

\noindent The diagonal maps are closed under precomposition with maps from
$\widehat{\mathbf{X}}$ and postcomposition with maps from $\widehat
{\mathbf{A}}$. Associativity follows from the associativity in the ambient
category $\mathbf{X}\times\mathbf{A}$.

The representation is accomplished essentially by putting a $\widehat{hat}$ on
objects and morphisms embedded in $\mathbf{X}\times\mathbf{A}$. The categories
$\mathbf{X}$ and $\mathbf{A}$ are represented respectively by the subcategory
$\widehat{\mathbf{X}}$ with objects $\widehat{x}=(x,Fx)$ and morphisms
$\widehat{f}=(f,Ff)$ and by the subcategory $\widehat{\mathbf{A}}$ with
objects $\widehat{a}=(Ga,a)$ and morphisms $\widehat{g}=(Gg,g)$. The
\textit{twist functor} $\left(  F,G\right)  :\mathbf{X}\times\mathbf{A}%
\rightarrow\mathbf{X}\times\mathbf{A}$ defined by $\left(  F,G\right)  \left(
x,a\right)  =\left(  Ga,Fx\right)  $ (and similarly for morphisms) restricted
to $\widehat{\mathbf{X}}\cong\mathbf{X}$ is $\widehat{F}$ which has the action
of $F$, i.e., $\widehat{F}\widehat{x}=(F,G)(x,Fx)=(GFx,Fx)=\widehat{Fx}%
\in\widehat{\mathbf{A}}$ and similarly for morphisms. The twist functor
restricted to $\widehat{\mathbf{A}}\cong\mathbf{A}$ yields $\widehat{G}$ which
has the action of $G$, i.e., $\widehat{G}\widehat{a}%
=(F,G)(Ga,a)=(Ga,FGa)=\widehat{Ga}\in\widehat{\mathbf{X}}$ and similarly for
morphisms. These functors provide representations on the left and right of the
abstract het-bifunctor $\operatorname*{Het}(\widehat{x},\widehat
{a})=\{\widehat{x}\overset{(f,f^{\ast})}{\longrightarrow}\widehat{a}\}$, i.e.,
the natural isomorphism

\begin{center}
$\operatorname*{Hom}_{\widehat{\mathbf{A}}}(\widehat{F}\widehat{x},\widehat
{a})\cong\operatorname*{Het}(\widehat{x},\widehat{a})\cong\operatorname*{Hom}%
_{\widehat{\mathbf{X}}}(\widehat{x},\widehat{G}\widehat{a})$.
\end{center}

\ This birepresentation of the abstract het-bifunctor gives an isomorphic copy
of the original adjunction between the isomorphic copies $\widehat{\mathbf{X}%
}$ and $\widehat{\mathbf{A}}$ of the original categories. This hom-pair
representation is summarized in the following:

\begin{description}
\item[\textbf{Adjunction Representation Theorem}:] Every adjunction
$F:\mathbf{X}\rightleftarrows\mathbf{A}:G$ can be represented (up to
isomorphism) as arising from the left and right representing universals of a
het-bifunctor $\operatorname*{Het}:\widehat{\mathbf{X}}^{op}\times
\widehat{\mathbf{A}}\rightarrow\mathbf{Set}$ giving the heteromorphisms from
the objects in a category $\widehat{\mathbf{X}}\cong\mathbf{X}$ to the objects
in a category $\widehat{\mathbf{A}}\cong\mathbf{A}$.\footnote{In a historical
note \cite[p. 103]{mac:cwm}, MacLane noted that Bourbaki \textquotedblleft
missed\textquotedblright\ the notion of an adjunction because Bourbaki focused
on the left representations of bifunctors $W:\mathbf{X}^{op}\times
\mathbf{A}\rightarrow\mathbf{Sets}$. MacLane remarks that given $G:\mathbf{A}%
\rightarrow\mathbf{X}$, they should have taken $W(x,a)=Hom_{\mathbf{X}}(x,Ga)$
and then focused on \textquotedblleft the symmetry of the adjunction
problem\textquotedblright\ to find $Fx$ so that $Hom_{\mathbf{A}}(Fx,a)\cong
Hom_{\mathbf{X}}(x,Ga)$. But MacLane thus missed the completely symmetrical
adjunction problem which is: given $W\left(  x,a\right)  $, find both $Ga$ and
$Fx$ such that $Hom_{\mathbf{A}}(Fx,a)\cong W(x,a)\cong Hom_{\mathbf{X}%
}(x,Ga)$.}
\end{description}

\subsection{Het Adjunctive Squares}

\noindent We previously used the representations of $\operatorname*{Het}(x,a)$
to pick out universal elements, the het unit $h_{x}\in\operatorname*{Het}%
(x,Fx)$ and the het counit $e_{a}\in\operatorname*{Het}(Ga,a)$, as the
respective adjoint correlates of $1_{Fx}$ and $1_{Ga}$ under the isomorphisms.
We showed that from the birepresentation of $\operatorname*{Het}(x,a)$, any
chimera morphism $x\overset{c}{\Rightarrow}a$ in $\operatorname*{Het}(x,a)$
would have two factorizations: $g(c)h_{x}=c=e_{a}f(c)$. This two
factorizations are spliced together along the main diagonal $c:x\Rightarrow a$
to form the het (commutative) adjunctive square.\pagebreak

\begin{center}
$%
\begin{array}
[c]{ccccc}
& x & \overset{f(c)}{\longrightarrow} & Ga & \\
h_{x} & \Downarrow & \searrow^{c} & \Downarrow & e_{a}\\
& Fx & \overset{g(c)}{\longrightarrow} & a &
\end{array}
\medskip$

Het Adjunctive Square\footnote{For typographical reasons, the diagonal
heteromorphism $c:x\Rightarrow a$ is represented as a single arrow rather than
a double arrow.}
\end{center}

Sometimes the two adjoint transposes are written vertically as in a
Gentzen-style rule of inference:

\begin{center}
\underline{$x\rightarrow Ga$}

$Fx\rightarrow a$

Gentzen-style presentation of an adjunction
\end{center}

\noindent This can be seen as a proto-het-adjunctive square without the
vertical morphisms---at least when the homomorphism involving the left adjoint
is on the bottom.

Some of the rigmarole of the conventional treatment of adjoints (\textit{sans}
chimeras) is only necessary because of the restriction to morphisms within one
category or the other. For instance, the UMP for the hom unit $\eta
_{x}:x\rightarrow GFx$ is that given any morphism $f:x\rightarrow Ga$ in
$\mathbf{X}$, there is a unique morphism $g=f^{\ast}:Fx\rightarrow a$ in the
other category $\mathbf{A}$ such that $G$-functorial image back in the
original category $\mathbf{X}$ gives the factorization of $f$ through the
unit: $x\overset{f}{\longrightarrow}Ga=x\overset{\eta_{x}}{\longrightarrow
}GFx\overset{Gf^{\ast}}{\longrightarrow}Ga$. The UMP has to go back and forth
between homomorphisms in the two categories because it avoids mention of the
heteromorphisms between the categories. The universal mapping property for the
het unit $h_{x}:x\Rightarrow Fx$ is much simpler (i.e., no $G$ and no over and
back). Given any heteromorphism $c:x\Rightarrow a$, there is a unique
homomorphism $g\left(  c\right)  :Fx\rightarrow a$ in the codomain category
$\mathbf{A}$ such that $x\overset{c}{\Rightarrow}a=x\overset{h_{x}%
}{\Rightarrow}Fx\overset{g\left(  c\right)  }{\longrightarrow}a$.

For instance, in the \textquotedblleft old days\textquotedblright\ (before
category theory), one might have stated the universal mapping property of the
free group $Fx$ on a set $x$ by saying that for any map $c:x\Rightarrow a$
from $x$ into a group $a$, there is a unique group homomorphism $g\left(
c\right)  :Fx\rightarrow a$ that preserves the action of $c$ on the generators
$x$, i.e., such that $x\overset{c}{\Rightarrow}a=x\hookrightarrow
Fx\overset{g\left(  c\right)  }{\longrightarrow}a$. That is just the left
half-adjunction part of the free-group adjunction. There is nothing sloppy or
`wrong' in that old way of stating the universal mapping property.

Dually for the hom counit, given any morphism $g:Fx\rightarrow a$ in
$\mathbf{A}$, there is a unique morphism $f=g^{\ast}:x\rightarrow Ga$ in the
other category $X$, such that the $F$-functorial image back in the original
category $A$ gives the factorization of $g$ though the counit: $Fx\overset
{g}{\longrightarrow}a=Fx\overset{Fg^{\ast}}{\longrightarrow}FGa\overset
{\varepsilon_{a}}{\longrightarrow}a$. For the het counit, given any
heteromorphism $c:x\Rightarrow a$, there is a unique homomorphism $f\left(
c\right)  :x\rightarrow Ga$ in the domain category $\mathbf{X}$ such that
$x\overset{c}{\Rightarrow}a=x\overset{f\left(  c\right)  }{\longrightarrow
}Ga\overset{e_{a}}{\Rightarrow}a$. Putting these two het UMPs together yields
the het adjunctive square diagram, just as previously putting the two hom UMPs
together yielded the hom-pair adjunctive square diagram.

\subsection{Het Natural Transformations}

One of the main motivations for category theory was to mathematically
characterize the intuitive notion of naturality for homomorphisms as in the
standard example of the canonical linear homomorphism embedded a vector space
into its double dual. Many heteromorphisms are rather arbitrary but certain
ones are quite canonical so we should be able to mathematically characterize
that canonicalness or naturality just as we do for homomorphisms. Indeed, the
notion of a natural transformation immediately generalizes to functors with
different codomains by taking the components to be heteromorphisms. \ Given
functors $F:\mathbf{X}\rightarrow\mathbf{A}$ and $H:\mathbf{X}\rightarrow
\mathbf{B}$ with a common domain \textit{and} given a het-bifunctor
$\operatorname*{Het}:\mathbf{A}^{op}\times\mathbf{B}\rightarrow\mathbf{Set}$,
a \textit{chimera }or\textit{ het natural transformation relative to
}$\operatorname*{Het}$, $\varphi:F\Rightarrow H$, is given by a set of
heteromorphisms $\{\varphi_{x}\in\operatorname*{Het}(Fx,Hx)\}$ indexed by the
objects of $\mathbf{X}$ such that for any $j:x\rightarrow x^{\prime}$ the
following diagram commutes.

\begin{center}
$%
\begin{array}
[c]{ccccc}
& Fx & \overset{\varphi_{x}}{\Longrightarrow} & Hx & \\
Fj & \downarrow &  & \downarrow & Hj\\
& Fx^{\prime} & \overset{\varphi_{x^{\prime}}}{\Longrightarrow} & Hx^{\prime}
&
\end{array}
\nolinebreak$

Het natural transformation
\end{center}

\noindent As with any commutative diagram involving heteromorphisms,
composition and commutativity are defined using the het-bifunctor (similar
remarks apply to any ordinary commutative hom diagram where it is the
hom-bifunctor behind the scenes). For instance, the above commutative squares
which define het natural transformations unpack as the following
behind-the-scenes commutative squares in $\mathbf{Set}$ for the underlying het-bifunctor.

\begin{center}
$%
\begin{array}
[c]{ccccc}
&  & \varphi_{x} &  & \\
& 1 & \longrightarrow & \operatorname*{Het}(Fx,Hx) & \\
\varphi_{x^{\prime}} & \downarrow &  & \downarrow & \operatorname*{Het}%
(Fx,Hj)\\
& \operatorname*{Het}(Fx^{\prime},Hx^{\prime}) & \longrightarrow &
\operatorname*{Het}(Fx,Hx^{\prime}) & \\
&  & \operatorname*{Het}(Fj,Hx^{\prime}) &  &
\end{array}
$
\end{center}

\noindent The composition $Fx\overset{\varphi_{x}}{\Longrightarrow}%
Hx\overset{Hj}{\longrightarrow}Hx^{\prime}$ is $\operatorname*{Het}%
(Fx,Hj)(\varphi_{x})\in\operatorname*{Het}(Fx,Hx^{\prime})$, the composition
$Fx\overset{Fj}{\longrightarrow}Fx^{\prime}\overset{\varphi_{x^{\prime}}%
}{\Longrightarrow}Hx^{\prime}$ is $\operatorname*{Het}(Fj,Hx^{\prime}%
)(\varphi_{x^{\prime}})\in\operatorname*{Het}(Fx,Hx^{\prime})$, and
commutativity means they are the same element of $\operatorname*{Het}%
(Fx,Hx^{\prime})$. These het natural transformations do not compose like the
morphisms in a functor category but they are acted upon by the natural
transformations in the functor categories on each side to yield het natural transformations.

There are het natural transformations each way between any functor and the
identity on its domain if the functor itself is used to define the appropriate
het-bifunctor. That is, given \textit{any} functor $F:\mathbf{X}%
\rightarrow\mathbf{A}$, there is a het natural transformation $1_{\mathbf{X}%
}\Rightarrow F$ relative to the bifunctor defined as $\operatorname*{Het}%
(x,a)=\operatorname*{Hom}_{\mathbf{A}}(Fx,a)$ as well as a het natural
transformation $F\Rightarrow1_{\mathbf{X}}$ relative to $\operatorname*{Het}%
(a,x)=\operatorname*{Hom}_{\mathbf{A}}(a,Fx)$.

Het natural transformations `in effect' already occur with reflective (or
coreflective) subcategories. A subcategory $\mathbf{A}$ of a category
$\mathbf{B}$ is a \textit{reflective subcategory} if the inclusion functor
$K:\mathbf{A}\hookrightarrow\mathbf{B}$ has a left adjoint. For any such
reflective adjunctions, the heteromorphisms $\operatorname*{Het}(b,a)$ are the
$\mathbf{B}$-morphisms with their heads in the subcategory $\mathbf{A}$ so the
representation on the right $\operatorname*{Het}(b,a)\cong\operatorname*{Hom}%
_{\mathbf{B}}(b,Ka)$ is trivial. The left adjoint $F:\mathbf{B}\rightarrow
\mathbf{A}$ gives the representation on the left: $\operatorname*{Hom}%
_{\mathbf{A}}(Fb,a)\cong\operatorname*{Het}(b,a)\cong\operatorname*{Hom}%
_{\mathbf{B}}(b,Ka)$. Then it is perfectly `natural' to see the unit of the
adjunction as defining a natural transformation $\eta:1_{\mathbf{B}%
}\Rightarrow F$ but that is actually a het natural transformation (since the
codomain of $F$ is $\mathbf{A}$). Hence the conventional (\textquotedblleft
heterophobic\textquotedblright?) treatment (e.g., \cite[p. 89]{mac:cwm}) is to
define another functor $R$ with the same domain and values on objects and
morphisms as $F$ except that its codomain is taken to be $\mathbf{B}$ so that
we can then have a hom natural transformation $\eta:1_{\mathbf{X}}\rightarrow
R$ between two functors with the same codomain. Similar remarks hold for the
dual coreflective case where the inclusion functor has a right adjoint and
where the heteromorphisms are turned around, i.e., are $\mathbf{B}$-morphisms
with their tail in the subcategory $\mathbf{A}$.

Given any adjunction isomorphism $\operatorname*{Hom}_{\mathbf{A}}%
(Fx,a)\cong\operatorname*{Het}(x,a)\cong\operatorname*{Hom}_{\mathbf{X}%
}(x,Ga)$, the adjoint correlates of the identities $1_{Fx}\in
\operatorname*{Hom}_{\mathbf{A}}\left(  Fx,Fx\right)  $ are the het units
$h_{x}\in\operatorname*{Het}\left(  x,Fx\right)  $ and the hom units $\eta
_{x}\in\operatorname*{Hom}_{\mathbf{X}}(x,GFx)$. The het units together give
the het natural transformation $h:1_{\mathbf{X}}\Rightarrow F$ while the hom
units give the hom natural transformation $\eta:1_{\mathbf{X}}\rightarrow GF$.
The adjoint correlates of the identities $1_{Ga}\in\operatorname*{Hom}%
_{\mathbf{X}}(Ga,Ga)$ are the het counits $e_{a}\in\operatorname*{Het}\left(
Ga,a\right)  $ and the hom counits $\varepsilon_{a}\in\operatorname*{Hom}%
_{\mathbf{A}}\left(  FGa,a\right)  .$ The het counits together give the het
natural transformation $e:G\Rightarrow1_{\mathbf{A}}$ while the hom counits
give the hom natural transformation $\varepsilon:FG\rightarrow1_{\mathbf{A}}$.

\section{Examples}

\subsection{The Product Adjunction}

Let $\mathbf{X}$ be the category $\mathbf{Set}$ of sets and let $\mathbf{A}$
be the category $\mathbf{Set}^{2}=\mathbf{Set}\times\mathbf{Set}$ of ordered
pairs of sets. A heteromorphism from a set to a pair of sets is a pair of set
maps with a common domain $\left(  f_{1},f_{2}\right)  :W\Rightarrow(X,Y)$
which is called a \textit{cone}. The het-bifunctor is given by
$\operatorname*{Het}\left(  W,\left(  X,Y\right)  \right)  =\left\{
W\Rightarrow\left(  X,Y\right)  \right\}  $, the set of all cones from $W$ to
$\left(  X,Y\right)  $. To construct a representation on the right, suppose we
are given a pair of sets $\left(  X,Y\right)  \in\mathbf{Set}^{2}$. How could
one construct a set, to be denoted $X\times Y$, such that all cones
$W\Rightarrow\left(  X,Y\right)  $ from any set $W$ could be represented by
set functions (morphisms within $\mathbf{Set}$) $W\rightarrow X\times Y$? In
the \textquotedblleft atomic\textquotedblright\ case of $W=1$ (the one element
set), a 1-cone $1\Rightarrow\left(  X,Y\right)  $ would just pick out an
ordered pair $\left(  x,y\right)  $ of elements, the first from $X$ and the
second from $Y$. Any cone $W\Rightarrow\left(  X,Y\right)  $ would just pick
out a set of pairs of elements. Hence the universal object would have to be
the set $\left\{  \left(  x,y\right)  :x\in X,y\in Y\right\}  $ of
\textit{all} such pairs which yields the Cartesian product of sets $X\times
Y$. The assignment of that set to each pair of sets gives the right adjoint
$G:\mathbf{Set}^{2}\rightarrow\mathbf{Set}$ where $G\left(  \left(
X,Y\right)  \right)  =X\times Y$ (and similarly for morphisms). The het counit
$e_{\left(  X,Y\right)  }:X\times Y\Rightarrow\left(  X,Y\right)  $
canonically takes each ordered pair $\left(  x,y\right)  $ as a single element
in $X\times Y$ to that pair of elements in $\left(  X,Y\right)  $. The
universal mapping property of the Cartesian product $X\times Y$ then holds;
given any set $W$ and a cone $\left(  f_{1},f_{2}\right)  :W\Rightarrow\left(
X,Y\right)  $, there is a unique set function $\left\langle f_{1}%
,f_{2}\right\rangle :W\rightarrow X\times Y$ defined by $\left\langle
f_{1},f_{2}\right\rangle \left(  w\right)  =\left(  f_{1}\left(  w\right)
,f_{2}\left(  w\right)  \right)  $ that factors the cone through het counit:

\begin{center}
$%
\begin{array}
[c]{cccc}%
W & \overset{\left\langle f_{1},f_{2}\right\rangle }{\longrightarrow} &
X\times Y & \\
\left(  f_{1},f_{2}\right)  & \searrow & \Downarrow & e_{\left(  X,Y\right)
}\\
&  & \left(  X,Y\right)  &
\end{array}
\nolinebreak$

Right half-adjunction of the product adjunction.
\end{center}

Fixing $W$ in the domain category, how could we find a universal object in
$\mathbf{Set}^{2}$ so that all heteromorphisms $\left(  f_{1},f_{2}\right)
:W\Rightarrow\left(  X,Y\right)  $ could be uniquely factored through it. The
obvious suggestion is the pair $\left(  W,W\right)  $ which defines a functor
$F:\mathbf{Set}\rightarrow\mathbf{Set}^{2}$ and where the het unit
$h_{W}:W\Rightarrow\left(  W,W\right)  $ is just the pair of identity maps
$h_{W}=(1_{W},1_{W})$. Then for each cone $\left(  f_{1},f_{2}\right)
:W\Rightarrow\left(  X,Y\right)  $, there is a unique pair of maps, also
denoted $\left(  f_{1},f_{2}\right)  :\left(  W,W\right)  \rightarrow\left(
X,Y\right)  $, which are a morphism in $\mathbf{Set}^{2}$ and which factors
the cone through the het unit:

\begin{center}
$%
\begin{array}
[c]{cccc}
& W &  & \\
h_{W} & \Downarrow & \searrow & \left(  f_{1},f_{2}\right) \\
& \left(  W,W\right)  & \overset{\left(  f_{1},f_{2}\right)  }{\longrightarrow
} & \left(  X,Y\right)
\end{array}
\nolinebreak$

Left half-adjunction of the product adjunction.
\end{center}

Splicing the two half-adjunctions along the diagonal gives the:

\begin{center}
$%
\begin{array}
[c]{ccccc}
& W & \overset{\left\langle f_{1},f_{2}\right\rangle }{\longrightarrow} &
X\times Y & \\
h_{W} & \Downarrow & \searrow^{\left(  f_{1},f_{2}\right)  } & \Downarrow &
e_{\left(  X,Y\right)  }\\
& \left(  W,W\right)  & \overset{\left(  f_{1},f_{2}\right)  }{\longrightarrow
} & \left(  X,Y\right)  &
\end{array}
\nolinebreak$

Het adjunctive square for the product adjunction.
\end{center}

The two factor maps on the top and bottom are uniquely associated with the
diagonal cones, and the isomorphism is natural so that we have natural
isomorphisms between the hom-bifunctors and the het-bifunctor:

\begin{center}
$\operatorname*{Hom}_{\mathbf{Set}^{2}}(\left(  W,W\right)  ,(X,Y))\cong%
\operatorname*{Het}\left(  W,\left(  X,Y\right)  \right)  \cong%
\operatorname*{Hom}_{\mathbf{Set}}\left(  W,X\times Y\right)  $.
\end{center}

\subsection{Limits in Sets}

Let $\mathbf{D}$ be a small (diagram) category and $D:\mathbf{D}{\rightarrow}$
$\mathbf{Set}$ a functor considered as a diagram in the category of
$\mathbf{Set}$. Limits in sets generalize the previous example where
$\mathbf{D}=2$. The diagram $D$ is in the functor category $\mathbf{Set}%
^{\mathbf{D}}$ where the morphisms are natural transformation between the
functors. A heteromorphism from a set $W$ to a diagram functor $D$ is
concretely given by a \textit{cone} $c:W\Rightarrow D$ which is defined as a
set of maps $\{W\overset{c_{i}}{\longrightarrow}D_{i}\}_{i\in Ob\left(
\mathbf{D}\right)  }$ (where $Ob\left(  \mathbf{D}\right)  $ is the set of
objects of $\mathbf{D}$) such that for any morphism $\alpha:i\rightarrow j$ in
$\mathbf{D}$, $W\overset{c_{i}}{\longrightarrow}D_{i}\overset{D_{\alpha}%
}{\longrightarrow}D_{j}=W\overset{c_{j}}{\longrightarrow}D_{j}$. The
adjunction is then given by the birepresentation of the het-bifunctor where
$\operatorname*{Het}(W,D)=\{W\Rightarrow D\}$ is the set of cones from the set
$W$ to the diagram functor $D$.

To construct the limit functor (right adjoint), take the product $\prod_{i\in
Ob\left(  \mathbf{D}\right)  }D_{i}$ and then take $LimD$ as the set of
elements $\left(  ...,x_{i},...\right)  $ of the product such that for any
morphism $\alpha:i\rightarrow j$ in $\mathbf{D}$, $D_{\alpha}(x_{i})=x_{j}$.

Since any cone $c:W\Rightarrow D=\left\{  W\overset{c_{i}}{\longrightarrow
}D_{i}\right\}  _{i\in Ob\left(  \mathbf{D}\right)  }$ would carry each
element $w\in W$ to a compatible set of elements from the $\left\{
D_{i}\right\}  $, the adjoint correlate map $f\left(  c\right)  :W\rightarrow
LimD$ would just carry $w$ to element $\left(  ...,c_{i}(w\right)  ,...)$ in
the set $LimD$. The het counit $e_{D}:LimD\Rightarrow D$ would be the cone
$\left\{  LimD\overset{\left(  e_{D}\right)  _{i}}{\longrightarrow}%
D_{i}\right\}  $ of maps such that $\left(  e_{D}\right)  _{i}:LimD\rightarrow
D_{i}$ carries each element $\left(  ...,x_{i},...\right)  $ of the set $LimD$
to the element $x_{i}\in D_{i}$ picked out in the $i^{th}$ component of
$\left(  ...,x_{i},...\right)  $. The universal mapping property for the het
counit $e_{D}$ is that given any cone $c:W\Rightarrow D$ (for any set $W$),
there is a unique set morphism $f\left(  c\right)  :W\rightarrow LimD$
(constructed above) such that $W\overset{f\left(  c\right)  }{\longrightarrow
}LimD\overset{e_{D}}{\Rightarrow}D=W\overset{c}{\Rightarrow}D$ which might be
diagrammed as:

\begin{center}
$%
\begin{array}
[c]{cccc}%
W & \overset{f(c)}{\longrightarrow} & LimD & \\
c & \searrow & \Downarrow & e_{D}\\
&  & D &
\end{array}
\nolinebreak$

Right half-adjunction of the limits adjunction.
\end{center}

Fixing $W$, how might we construct a functor $\mathbf{D}\rightarrow
\mathbf{Set}$ such that all cones $c:W\Rightarrow D$ might be factored through
it? Since the only given data is $W$, the obvious thing to try is the constant
or diagonal functor $\Delta W:\mathbf{D}\rightarrow\mathbf{Set}$ which takes
each object $i$ in $\mathbf{D}$ to $W$ and each morphism $\alpha:i\rightarrow
j$ to the identity map $1_{W}$. Given a cone $c:W\Rightarrow D$, the obvious
natural transformation $g\left(  c\right)  $ from the diagonal functor $\Delta
W$ to $D$ is given by the same set of maps $g\left(  c\right)  =\left\{
W\overset{c_{i}}{\longrightarrow}D_{i}\right\}  _{i\in Ob\left(
\mathbf{D}\right)  }$ which constitutes a morphism in $\mathbf{Set}%
^{\mathbf{D}}$.\footnote{It should be carefully noted that much of the
literature refers to this natural transformation morphism within the category
$\mathbf{Set}^{\mathbf{D}}$ interchangably with the cone $W\Rightarrow D$
which is a heteromorphism from a set to a functor.} The het unit
$h_{W}:W\Rightarrow\Delta W$ is the set of identity maps $\left\{
W\overset{\left(  h_{W}\right)  _{i}=1_{W}}{\longrightarrow}\left(  \Delta
W\right)  _{i}=W\right\}  $. The (rather trivial) universal mapping property
of the het unit $h_{W}$ is that given any cone $c:W\Rightarrow D$, there
exists a unique natural transformation $g\left(  c\right)  :\Delta
W\rightarrow D$ which factors the cone through the het unit:

\begin{center}
$%
\begin{array}
[c]{cccc}
& W &  & \\
h_{W} & \Downarrow & \searrow & c\\
& \Delta W & \overset{g(c)}{\longrightarrow} & D
\end{array}
\nolinebreak$

Left half-adjunction of the limit adjunction.
\end{center}

\noindent Splicing the two half-adjunctions together along the main diagonal
gives the:

\begin{center}
$%
\begin{array}
[c]{ccccc}
& W & \overset{f(c)}{\longrightarrow} & LimD & \\
h_{W} & \Downarrow & \searrow^{c} & \Downarrow & e_{D}\\
& \Delta W & \overset{g(c)}{\longrightarrow} & D &
\end{array}
\nolinebreak$

Het adjunctive square for the limit adjunction.
\end{center}

The two factor maps on the top and bottom are uniquely associated with the
diagonal cones, and the isomorphism is natural so that we have natural
isomorphisms between the hom-bifunctors and the het-bifunctor:

\begin{center}
$\operatorname*{Hom}_{\mathbf{Set}^{D}}(\Delta W,D)\cong\operatorname*{Het}%
\left(  W,D\right)  \cong\operatorname*{Hom}_{\mathbf{Set}}\left(
W,LimD\right)  $.
\end{center}

\subsection{Colimits in Sets}

The duals to products are coproducts which in $Set$ is the disjoint union of
sets. Instead of developing that special case, we move directly to the
adjunction for colimits in $\mathbf{Set}$ which is dual to the limits
adjunction. The argument here (as with limits) would work for any other
cocomplete (or complete) category of algebras replacing the category of sets.
Given the same data as in the previous section, the diagonal functor
$\Delta:\mathbf{Set}\rightarrow\mathbf{Set}^{\mathbf{D}}$ is defined as before
and it has a \textit{left} adjoint $\operatorname*{Colim}:\mathbf{Set}%
^{\mathbf{D}}\rightarrow\mathbf{Set}$.

For this adjunction, a heteromorphism from a diagram functor $D$ to a set $Z$
is a \textit{cocone} $c:D\Rightarrow Z$ which is defined as a set of maps
$\{D_{i}\overset{c_{i}}{\longrightarrow}Z\}_{i\in Ob\left(  \mathbf{D}\right)
}$ such that for any morphism $\alpha:i\rightarrow j$ in $\mathbf{D}$,
$D_{i}\overset{D_{\alpha}}{\longrightarrow}D_{j}\overset{c_{j}}%
{\longrightarrow}Z=D_{i}\overset{c_{i}}{\longrightarrow}Z$. The adjunction is
then given by the birepresentations of the het-bifunctor where
$\operatorname*{Het}(D,Z)=\{D\Rightarrow Z\}$ is the set of cocones from the
diagram functor $D$ to the set $Z$.

Fixing a diagram functor $D:\mathbf{D}\rightarrow Set$ in $\mathbf{Set}%
^{\mathbf{D}}$, how can we construct a universal object $\operatorname*{Colim}%
D$ so that any cocone $c:D\Rightarrow Z$ can be factored through it? Instead
of taking the product $\prod_{i\in Ob\left(  \mathbf{D}\right)  }D_{i}$, take
the dual construction of the coproduct which is the disjoint union
$\coprod\limits_{i\in Ob\left(  \mathbf{D}\right)  }D_{i}$. Then instead of
the subset of compatible elements, take $\operatorname*{Colim}D$ as the
quotient set by the compatibility equivalence relation $x_{i}\sim x_{j}$ if
$D_{\alpha}(x_{i})=x_{j}$ for any morphism $D_{\alpha}$ between the $D_{i}$s.

To construct the het unit cocone $h_{D}:D\Rightarrow\operatorname*{Colim}D$,
define each map $\left(  h_{D}\right)  _{i}:D_{i}\rightarrow
\operatorname*{Colim}D$ by taking each element $x\in D_{i}$ to its equivalence
class in $\operatorname*{Colim}D$. All the maps in $h_{D}=\left\{
D_{i}\overset{\left(  h_{D}\right)  _{i}}{\longrightarrow}%
\operatorname*{Colim}D\right\}  _{i\in Ob\left(  \mathbf{D}\right)  }$ will
commute with the maps $D\alpha:D_{i}\rightarrow D_{j}$ for $\alpha
:i\rightarrow j$ in $\mathbf{D}$ by the construction of the coelements so that
$h_{D}$ is a cocone $D\Rightarrow\operatorname*{Colim}D$. Given a cocone
$c:D\Rightarrow Z=\left\{  D_{i}\overset{c_{i}}{\longrightarrow}Z\right\}
_{i\in Ob\left(  \mathbf{D}\right)  }$, an $x\in D_{i}$ must be carried by
$c_{i}$ to the same element $z$ of $Z$ as all the other elements $x^{\prime
}\in D_{j}$ in the equivalence class of $x$ are carried by their maps
$c_{j}:D_{j}\rightarrow Z$, so the factor set map $g\left(  c\right)
:\operatorname*{Colim}D\rightarrow Z$ would just carry that equivalence class
of $x$ as an element of $\operatorname*{Colim}D$ to that $z\in Z$. This factor
map uniquely completes the following diagram:

\begin{center}
$%
\begin{array}
[c]{cccc}
& D &  & \\
h_{D} & \Downarrow & \searrow & c\\
& \operatorname*{Colim}D & \overset{g\left(  c\right)  }{\longrightarrow} & Z
\end{array}
\nolinebreak$

Left half-adjunction for colimit adjunction.
\end{center}

Fixing $Z$, the same diagonal functor $\Delta:\mathbf{Set}\rightarrow
\mathbf{Set}^{\mathbf{D}}$ gives a universal object $\Delta Z$ in
$\mathbf{Set}^{\mathbf{D}}$. The het counit $e_{Z}$ is the cocone
$e_{Z}=\left\{  \left(  \Delta Z\right)  _{i}=Z\overset{1_{Z}}{\rightarrow
}Z\right\}  _{i\in Ob\left(  \mathbf{D}\right)  }$ of identity maps. Given any
cocone $c:D\Rightarrow Z=\left\{  D_{i}\overset{c_{i}}{\rightarrow}Z\right\}
_{i\in Ob\left(  \mathbf{D}\right)  }$, the adjoint correlate natural
transformation $f\left(  c\right)  :D\rightarrow\Delta Z$ would be defined by
the same set of maps $\left\{  D_{i}\overset{c_{i}}{\rightarrow}Z\right\}  $.
That is clearly the unique natural transformation $D\rightarrow\Delta Z$ to
make the following diagram commute:

\begin{center}
$%
\begin{array}
[c]{cccc}%
D & \overset{f\left(  c\right)  }{\longrightarrow} & \Delta Z & \\
c & \searrow & \Downarrow & e_{Z}\\
&  & Z &
\end{array}
\nolinebreak$

Right half-adjunction for colimit adjunction.
\end{center}

\noindent Splicing the two half-adjunctions together along the main diagonal
yields the:

\begin{center}
$%
\begin{array}
[c]{ccccc}
& D & \overset{f\left(  c\right)  }{\longrightarrow} & \Delta Z & \\
h_{D} & \Downarrow & \searrow^{c} & \Downarrow & e_{Z}\\
& \operatorname*{Colim}D & \overset{g\left(  c\right)  }{\longrightarrow} &
Z &
\end{array}
\nolinebreak$

Het adjunctive square for the colimit adjunction.
\end{center}

The two factor maps on the top and bottom are uniquely associated with the
diagonal cocone, and the isomorphism is natural so that we have natural
isomorphisms between the hom-bifunctors and the het-bifunctor:

\begin{center}
$\operatorname*{Hom}_{\mathbf{Set}}(\operatorname*{Colim}D,Z)\cong%
\operatorname*{Het}\left(  D,Z\right)  \cong\operatorname*{Hom}_{\mathbf{Set}%
^{\mathbf{D}}}\left(  D,\Delta Z\right)  $.
\end{center}

\subsection{Adjoints to Forgetful Functors}

Perhaps the most accessible adjunctions are the free-forgetful adjunctions
between $\mathbf{X}=\mathbf{Set}$ and a category of algebras such as the
category of groups $\mathbf{A}=\mathbf{Grps}$. \text{The right adjoint
}$G:\mathbf{A}\rightarrow\mathbf{X}$ forgets the group structure to give the
underlying set $GA$ of a group $A$. \ The left adjoint $F:\mathbf{X}%
\rightarrow\mathbf{A}$ gives the free group $FX$ generated by a set $X$.

For this adjunction, the heteromorphisms are any set functions $X\overset
{c}{\Rightarrow}A$ (with the codomain being a group $A$) and the het-bifunctor
is given by such functions: $\operatorname*{Het}(X,A)=\{X\Rightarrow A\}$
(with the obvious morphisms). A heteromorphism $c:X\Rightarrow A$ determines a
set map $f(c):X\rightarrow GA$ trivially and it determines a group
homomorphism $g(c):FX\rightarrow A$ by mapping the generators $x\in X$ to
their images $c\left(  x\right)  \in A$ and then mapping the other elements of
$FX$ as they must be mapped in order for $g\left(  c\right)  $ to be a group
homomorphism. The het unit $h_{X}:X\Rightarrow FX$ is insertion of the
generators into the free group and the het counit $e_{A}:GA\Rightarrow A$ is
just the retracting of the elements of the underlying set back to the group.
These factor maps $f\left(  c\right)  $ and $g\left(  c\right)  $ uniquely
complete the usual two half-adjunction triangles which together give the:

\begin{center}
$%
\begin{array}
[c]{ccccc}
& X & \overset{f\left(  c\right)  }{\longrightarrow} & GA & \\
h_{X} & \Downarrow & \searrow^{c} & \Downarrow & e_{A}\\
& FX & \overset{g\left(  c\right)  }{\longrightarrow} & A &
\end{array}
\nolinebreak$

Het adjunctive square for the free group adjunction.
\end{center}

These associations also give us the two representations:

\begin{center}
$\operatorname*{Hom}(FX,A)\cong\operatorname*{Het}(X,A)\cong%
\operatorname*{Hom}(X,GA)$.
\end{center}

In general, the existence of a left adjoint to $U:\mathbf{A}\rightarrow
\mathbf{Set}$ (i.e., a left representation of $\operatorname*{Het}%
(X,A)=\{X\Rightarrow A\}$) will depend on whether or not there is an
$\mathbf{A}$-object $FX$ with the least or minimal structure so that every
chimera $X\overset{c}{\Rightarrow}A$ will determine a unique representing
$\mathbf{A}$\textbf{-}morphism $g\left(  c\right)  :FX\rightarrow A$.

The existence of a \textit{right} adjoint to $U$ will depend on whether or not
for any set $X$ there is an $\mathbf{A}$-object $IX$ with the greatest or
maximum structure so that any chimera $A\Rightarrow X$ can be represented by
an $\mathbf{A}$-morphism $A\rightarrow IX$.

Consider the underlying set functor $U:\mathbf{Pos}\rightarrow\mathbf{Set}$
from the category of partially ordered sets (an ordering that is reflexive,
transitive, and anti-symmetric) with order-preserving maps to the category of
sets. It has a left adjoint since each set has a least partial order on it,
namely the discrete ordering. \ Hence any chimera function $X\overset
{c}{\Rightarrow}A$\ from a set $X$ to a partially ordered set or poset $A$
could be represented as a set function $X\overset{f\left(  c\right)
}{\longrightarrow}UA$ or as an order-preserving function $DX\overset{g\left(
c\right)  }{\longrightarrow}A$ where $DX$ gives the discrete ordering on $X$.
The functor giving the discrete partial ordering on a set is left adjoint to
the underlying set function.

In the other direction, one could take as a chimera any function $A\overset
{c}{\Rightarrow}X$ (from a poset $A$ to a set $X$)\ and it is represented on
the left by the ordinary set function $UA\overset{g\left(  c\right)
}{\longrightarrow}X$ so the left half-adjunction trivially exists:

\begin{center}
$%
\begin{array}
[c]{cccc}
& A & \overset{?}{\longrightarrow} & IX?\\
h_{A} & \Downarrow & \searrow^{c} & \downarrow?\\
& UA & \overset{g\left(  c\right)  }{\longrightarrow} & X
\end{array}
\nolinebreak$

Left half-adjunction (with no right half-adjunction).
\end{center}

But the underlying set functor $U$ does not have a right adjoint since there
is no maximal partial order $IX$ on $X$ so that any chimera $A\overset
{c}{\Rightarrow}X$ could be represented as an order-preserving function
$f\left(  c\right)  :A\rightarrow IX$. To receive all the possible orderings,
the ordering relation would have to go both ways between any two points which
would then be identified by the anti-symmetry condition so that $IX$ would
collapse to a single point and the factorization of $c$ through $IX$ would
fail.\footnote{Thanks to Vaughn Pratt for the example.} Thus poset-to-set
chimera $A\Rightarrow X$\ can only be represented on the left.

Relaxing the anti-symmetry condition, let $U:\mathbf{Ord}\rightarrow
\mathbf{Set}$ be the underlying set functor from the category of preordered
sets (reflexive and transitive orderings) to the category of sets. \ The
discrete ordering again gives a left adjoint. \ But now there is also a
maximal ordering on a set $X$, namely the `indiscrete' ordering $IX$ on $X$
(the `indiscriminate' or `chaotic' preorder on $X$) which has the ordering
relation both ways between any two points. \ Then a preorder-to-set chimera
morphism $A\Rightarrow X$\ (just a set function ignoring the ordering) can be
represented on the left as a set function $UA\overset{g\left(  c\right)
}{\longrightarrow}X$ and on the right as an order-preserving function
$A\overset{f\left(  c\right)  }{\longrightarrow}IX$ so that $U$ also has a
right adjoint $I$ and we have the following:

\begin{center}
$%
\begin{array}
[c]{ccccc}
& A & \overset{f\left(  c\right)  }{\longrightarrow} & IX & \\
h_{A} & \Downarrow & \searrow^{c} & \Downarrow & e_{X}\\
& UA & \overset{g\left(  c\right)  }{\longrightarrow} & X &
\end{array}
\nolinebreak$

Het adjunctive square for the indiscrete-underlying adjunction on preorders.
\end{center}

\subsection{Reflective Subcategories}

Suppose that $\mathbf{A}$ is a subcategory of $\mathbf{X}$ with $G:\mathbf{A}%
\hookrightarrow\mathbf{X}$ the inclusion functor and suppose that it has a
left adjoint $F:\mathbf{X}\rightarrow\mathbf{A}$. Then $\mathbf{A}$ is said to
be a \textit{reflective subcategory} of $\mathbf{X}$, the left adjoint $F$ is
the \textbf{reflector}, and the adjunction is called a \textit{reflection}:
$\operatorname*{Hom}_{\mathbf{A}}(Fx,a)\cong\operatorname*{Hom}_{\mathbf{X}%
}(x,Ga)$. For all reflections, the chimera morphisms are the morphisms
$x\Rightarrow a$ in the ambient category $\mathbf{X}$ with their heads in the
reflective subcategory $\mathbf{A}$. Hence the het-bifunctor would be:

\begin{center}
$\operatorname*{Het}(x,a)=\operatorname*{Hom}_{\mathbf{XA}}(x,a)$
\end{center}

\noindent where the $\mathbf{XA}$ subscript indicates that $x$ can be any
object in $\mathbf{X}$ but that $a$\ is any element of the subcategory
$\mathbf{A}$. Note the two ways of seeing any $c\in\operatorname*{Het}%
(x,a)=\operatorname*{Hom}_{\mathbf{XA}}(x,a)$. From one viewpoint,
$c\in\operatorname*{Hom}_{\mathbf{XA}}(x,a)\subseteq\operatorname*{Hom}%
_{\mathbf{X}}\left(  x,a\right)  $ so that $c$ is just a morphism inside the
category $\mathbf{X}$, but we also view it as a chimera with its tail in
$\mathbf{X}$ and head in $\mathbf{A}$. Since $G$ is the inclusion functor, it
just takes $a$\ as an element of $\mathbf{A}$ to itself as an element of
$\mathbf{X}$ and similarly for morphisms. Thus we insert $\operatorname*{Het}%
\left(  x,a\right)  $ in the middle to get the two representation isomorphisms:

\begin{center}
$\operatorname*{Hom}_{\mathbf{A}}(Fx,a)\cong\operatorname*{Het}(x,a)\cong%
\operatorname*{Hom}_{\mathbf{X}}(x,Ga)$.
\end{center}

There is also the dual case of a \textit{coreflective subcategory} where the
inclusion functor has a right adjoint and where the chimera morphisms are
turned around (tail in subcategory and head in ambient category). This case
will be used in the next section but here I will focus on reflective subcategories.

For an interesting example of a reflector dating back five centuries, we use
the modern mathematical formulation of double-entry bookkeeping
\cite{ell:mdeb} \cite{ell:dema}. Let $\mathbf{Ab}$ be the category of abelian
(i.e., commutative) groups where the operation is written as addition. Thus
$0$ is the identity element, $a+a^{\prime}=a^{\prime}+a$, and for each element
$a$, there is an element $-a$ such that $a+(-a)=0$. Let $\mathbf{CMon}$ be the
category of commutative monoids so the addition operation has the identity $0$
but does not necessarily have an inverse. Let $G:\mathbf{Ab}\hookrightarrow
\mathbf{CMon}$ be the inclusion functor.

In 1494, the mathematician Luca Pacioli published an accounting technique that
had been developed in practice during the 1400s and which became known as
\textit{double-entry bookkeeping} \cite{pac:deb}. In essence, the idea was to
do additive arithmetic with additive inverses using ordered pairs
$[x\ //\ x^{\prime}]$ of non-negative numbers called \textit{T-accounts}%
.\footnote{The double-slash separator was suggested by Pacioli.
{\textquotedblleft}At the beginning of each entry, we always provide `per',
because, first, the debtor must be given, and immediately after the creditor,
the one separated from the other by two little slanting parallels
(virgolette), thus, //,\dots.{\textquotedblright\ \cite[p. 43]{pac:deb}}} The
number on the left side was called the \textit{debit entry} and the number on
the right the \textit{credit entry}. T-accounts added by adding the
corresponding entries: $[x\ //\ x^{\prime}]+[y\ //\ y^{\prime}%
]=[x+y\ //\ x^{\prime}+y^{\prime}]$. Two T-accounts were deemed equal if their
cross-sums were equal (the additive version of the equal cross-multiples used
to define equality of multiplicative ordered pairs or fractions). Thus

\begin{center}
$[x\ //\ x^{\prime}]=[y\ //\ y^{\prime}]$ if $x+y^{\prime}=x^{\prime}+y$.
\end{center}

Hence the additive inverse was obtained by \textquotedblleft reversing the
entries\textquotedblright\ (as accountants say):

\begin{center}
$[x\ //\ x^{\prime}]+[x^{\prime}\ //\ x]=[x+x^{\prime}\ //\ x^{\prime
}+x]=[0\ //\ 0]$.
\end{center}

To obtain the reflector or left adjoint $F:\mathbf{CMon}\rightarrow
\mathbf{Ab}$ to $G$, we need only note that Pacioli was implicitly using the
fact that the normal addition of numbers is \textit{cancellative} in the sense
that $x+z=y+z$ implies $x=y$. Since commutative monoids do not in general have
that property we need only to tweak the definition of equality of T-accounts
\cite[p. 17]{bour:alg1}:

\begin{center}
$[x\ //\ x^{\prime}]=[y\ //\ y^{\prime}]$ if there is a $z$ such that
$x+y^{\prime}+z=x^{\prime}+y+z$.
\end{center}

This construction with the induced maps then yields a functor $F:\mathbf{CMon}%
\rightarrow\mathbf{Ab}$ that takes a commutative monoid $m$ to a commutative
group $Fm=P(m)$. The group $P(m)$ is usually called the \textquotedblleft
group of differences\textquotedblright\ or \textquotedblleft
inverse-completion\textquotedblright\ and, in algebraic geometry, its
generalization is called the \textquotedblleft Grothendieck
group.\textquotedblright\ However, due to about a half-millennium of priority,
we will call the additive group of differences the \textit{Pacioli group }of
the commutative monoid\textit{ }$m$. For any such $m$, the het unit
$h_{m}:m\Rightarrow Fm=P(m)$ which takes an element $x$ to the T-account
$[0\ //\ x]$ with that credit balance (the debit balance mapping would do just
as well).

For this adjunction, a heteromorphism $c:m\Rightarrow a$, is any
\textit{monoid homomorphism} from a commutative monoid $m$ to any abelian
group $a$ (being only a monoid homomorphism, it does not need to preserve any
inverses that might exist in $m$). The Pacioli group has the following
universality property: for any heteromorphism $c:m\Rightarrow a$, there is a
unique group homomorphism $g(c):Fm\rightarrow a$ such that $m\overset{h_{m}%
}{\Rightarrow}Fm\overset{g(c)}{\rightarrow}a=m\overset{c}{\Rightarrow}a$. The
group homomorphism factor map is: $g(c)([x^{\prime}%
\ //\ x])=c(x)+(-c(x^{\prime}))$. This establishes the other representation
isomorphism of the adjunction:

\begin{center}
$\operatorname*{Hom}_{\mathbf{Ab}}(Fm,a)\cong\operatorname*{Het}%
(m,a)\cong\operatorname*{Hom}_{\mathbf{CMon}}(m,Ga)$.
\end{center}

\subsection{The Special Case of Endo-Adjunctions}

Since the heteromorphic theory of adjoints is based on representing the
heteromorphisms between the objects of two different categories with
homomorphisms within each category, the case of an endo-adjunction all within
one category is clearly going to require some special attention. The
product-exponential adjunction in $\mathbf{Set}$ is an important example of
adjoint endo-functors $F:\mathbf{Set}\rightleftarrows\mathbf{Set}:G$. For any
fixed (non-empty) \textquotedblleft index\textquotedblright\ set $A$, the
product functor $F(-)=-\times A:\mathbf{Set}\rightarrow\mathbf{Set}$ has a
right adjoint $G(-)=(-)^{A}:\mathbf{Set}\rightarrow\mathbf{Set}$ which makes
$\mathbf{Set}$ a Cartesian-closed category. For any sets $X$ and $Y$, the
adjunction has the form: $\operatorname*{Hom}(X\times A,Y)\cong%
\operatorname*{Hom}(X,Y^{A})$.

Since both functors are endo-functors on $\mathbf{Set}$, we don't have the two
categories between which to have heteromorphisms. Moreover, we don't have the
expected canonical maps as the het unit or counit. For instance, the het unit
should be a canonical morphism $h_{X}:X\Rightarrow FX$ but if $FX=X\times A$,
there is no \textit{canonical} (het or otherwise) map $X\rightarrow X\times A$
(except in the special case where $A$ is a singleton). Similarly, the het
counit should be a canonical map $e_{Y}:GY\Rightarrow Y$ but if $GY=Y^{A}$
then there is no canonical (het or otherwise) map $Y^{A}\rightarrow Y$ (unless
$A$ is a singleton). Hence a special treatment is required. It consists of
showing that the endo-adjunction can be parsed in two ways as adjunctions each
of which is between different categories and then the heteromorphic theory applies.

The key is the following special case of a result by Freyd \cite[p.
83]{Freyd:ac}. Consider any endo-adjunction $F:\mathbf{C}\rightleftarrows
\mathbf{C}:G$ on a category $\mathbf{C}$ where the functors are assumed
one-one on objects. Then the image of $G$ is a subcategory of $\mathbf{C}$,
i.e., $\mathbf{\operatorname{Im}}\left(  G\right)  \hookrightarrow\mathbf{C}$,
and similarly the image of $F$ is also a subcategory of $\mathbf{C}$, i.e.,
$\mathbf{\operatorname{Im}}\left(  F\right)  \hookrightarrow\mathbf{C}$. The
operation of taking the $G$-image of the hom-set $\operatorname*{Hom}%
_{\mathbf{C}}(Fx,a)$ to obtain $\operatorname*{Hom}_{\mathbf{\operatorname{Im}%
}\left(  G\right)  }(GFx,Ga)$ is onto by construction. It is one-one on
objects by assumption and one-one on maps since if for $g,g^{\prime
}:Fx\rightarrow a$ and $Gg=Gg^{\prime}$, then $g=g^{\prime}$ by the uniqueness
of the factor map $f^{\ast}:Fx\rightarrow a$ to factor $x\overset
{f}{\rightarrow}Ga:=x\overset{\eta_{x}}{\longrightarrow}GFx\overset
{Gg=Gg^{\prime}}{\longrightarrow}Ga$ through the hom unit $\eta_{x}$. The
isomorphism is also natural in $x$ and $a$ so we have:

\begin{center}
$\operatorname*{Hom}_{\mathbf{\operatorname{Im}}\left(  G\right)
}(GFx,Ga)\cong\operatorname*{Hom}_{\mathbf{C}}(Fx,a)$.
\end{center}

If $x^{\prime}=Ga$, then using that isomorphism and the adjunction, we have:

\begin{center}
$\operatorname*{Hom}_{\mathbf{\operatorname{Im}}\left(  G\right)
}(GFx,x^{\prime})\cong\operatorname*{Hom}_{\mathbf{\operatorname{Im}}\left(
G\right)  }(GFx,Ga)\cong\operatorname*{Hom}_{\mathbf{C}}(Fx,a)\cong%
\operatorname*{Hom}_{\mathbf{C}}(x,Ga)\cong\operatorname*{Hom}_{\mathbf{C}%
}(x,x^{\prime})$.
\end{center}

\noindent Thus $\operatorname*{Hom}_{\mathbf{\operatorname{Im}}\left(
G\right)  }(GFx,x^{\prime})\cong$ $\operatorname*{Hom}_{\mathbf{C}%
}(x,x^{\prime})$ so that $GF:\mathbf{C}\rightarrow\mathbf{\operatorname{Im}%
(}G\mathbf{)}$ (where $G$ is construed as having the codomain
$\operatorname{Im}\left(  G\right)  $) is left adjoint to the inclusion
$\mathbf{\operatorname{Im}}\left(  G\right)  \hookrightarrow\mathbf{C}$ and
thus $\mathbf{\operatorname{Im}}\left(  G\right)  $ is a reflective
subcategory of $\mathbf{C}$.

Dually, we also have the natural isomorphism

\begin{center}
$\operatorname*{Hom}_{\mathbf{C}}(x,Ga)\cong\operatorname*{Hom}%
_{\mathbf{\operatorname{Im}}\left(  F\right)  }(Fx,FGa)$
\end{center}

\noindent by taking the $F$-image of $\operatorname*{Hom}_{\mathbf{C}}(x,Ga)$
and using the universal mapping property of the hom counit. If $a^{\prime}=Fx$
then using that isomorphism and the adjunction, we have:

\begin{center}
$\operatorname*{Hom}_{\mathbf{C}}(a^{\prime},a)\cong\operatorname*{Hom}%
_{\mathbf{C}}(Fx,a)\cong\operatorname*{Hom}_{\mathbf{C}}(x,Ga)\cong%
\operatorname*{Hom}_{\mathbf{\operatorname{Im}}\left(  F\right)
}(Fx,FGa)\cong\operatorname*{Hom}_{\mathbf{\operatorname{Im}}\left(  F\right)
}(a^{\prime},FGa)$.
\end{center}

\noindent Thus $\operatorname*{Hom}_{\mathbf{C}}(a^{\prime},a)\cong%
\operatorname*{Hom}_{\mathbf{\operatorname{Im}}\left(  F\right)  }(a^{\prime
},FGa)$ so that $FG:\mathbf{C}\rightarrow\mathbf{\operatorname{Im}}\left(
F\right)  $ (where $F$ is construed as having the codomain $\operatorname{Im}%
\left(  F\right)  $) is right adjoint to the inclusion
$\mathbf{\operatorname{Im}}\left(  F\right)  \hookrightarrow\mathbf{C}$ and
thus $\mathbf{\operatorname{Im}(}F\mathbf{)}$ is a coreflective subcategory of
$\mathbf{C}$.

Therefore the endo-adjunction can be analyzed or parsed into two adjunctions
between \textit{different} categories, a reflection and a coreflection. It was
previously noted that in the case of a reflection, i.e., a left adjoint to the
inclusion functor, heteromorphisms can be found as the morphisms with their
tails in the ambient category and their heads in the subcategory. For a
coreflection (right adjoint to the inclusion functor), the heteromorphisms
would be turned around, i.e., would have their tails in the subcategory and
their heads in the ambient category. The heteromorphic theory applies to each
of these adjunctions.

In the case of the exponential endo-adjunction on $\mathbf{Set}$,
$\mathbf{X}=\mathbf{Set}=\mathbf{A}$. To parse the adjunction as a reflection,
let $\mathbf{APower}$ be the subcategory of $G(-)=(-)^{A}$ images ($G$ is
one-one since $A$ is non-empty) so that $\mathbf{APower}\hookrightarrow
\mathbf{Set}$, and that inclusion functor has a left adjoint $GF(-)=(-\times
A)^{A}:\mathbf{Set}\rightarrow\mathbf{APower}$. Then the heteromorphisms are
those with their tail in $\mathbf{Set}$ and head in $\mathbf{APower}$, i.e.,
the morphisms of the form $X\rightarrow Y^{A}$. But now we have the het unit
and counit in accordance with the heteromorphic theory. The het unit
$h_{X}:X\Rightarrow GFX=(X\times A)^{A}$ is the canonical map that takes an
$x$ in $X$ to the function $(x,-):A\rightarrow X\times A$ which takes $a$ in
$A$ to $(x,a)\in X\times A$. This is the `same' as the ordinary unit $\eta
_{X}:X\rightarrow(X\times A)^{A}$ in the original product-exponential
adjunction, i.e., it is the `same' modulo the fact that $h_{x}$ is viewed as a
heteromorphism with its tail in $\mathbf{Set}$ and its head in
$\mathbf{APower}$ while the same set-map $\eta_{X}$ is viewed as a
homomorphism in $\mathbf{Set}$. Since the right adjoint in the reflective case
is the inclusion, the het counit $e_{Y^{A}}:Y^{A}\rightarrow Y^{A}$ is the
identity but seen as a heteromorphism from an object in $\mathbf{Set}$ to the
same object in $\mathbf{APower}$. The het adjunctive square then is the
following commutative diagram.

\begin{center}
$%
\begin{array}
[c]{ccccc}
& X & \overset{f}{\rightarrow} & Y^{A} & \\
h_{X} & \Downarrow & \searrow^{f} & \Downarrow & e_{Y^{A}}\\
& (X\times A)^{A} & \overset{(f^{\ast})^{A}}{\longrightarrow} & Y^{A} &
\end{array}
\nolinebreak$

Het adjunctive square for exponential endo-adjunction parsed as a reflection
\end{center}

As a coreflection, let $\mathbf{AProd}$ be the subcategory of $F(-)=-\times A$
images ($F$ is one-one since $A$ is non-empty) so that $\mathbf{AProd}%
\hookrightarrow\mathbf{Set}$, and that inclusion functor has a right adjoint
$FG(-)=(-)^{A}\times A:\mathbf{Set}\rightarrow\mathbf{AProd}$. Then the
heteromorphisms are those with their tail in $\mathbf{AProd}$ and their head
in $\mathbf{Set}$, i.e., the morphisms of the form $X\times A\rightarrow Y$.
And now we again have the het unit and counit. The het counit $e_{Y}%
:Y^{A}\times A=FGY\Rightarrow Y$ is the evaluation map which is the `same' as
the ordinary counit $\varepsilon_{Y}:Y^{A}\times A\rightarrow Y$ in the
original product-exponential adjunction. Since the left adjoint in the
coreflective case is the inclusion, the het unit $h_{X\times A}:X\times
A\Rightarrow X\times A$ is the identity (again seen as a heteromorphism). The
het adjunctive square is then the following commutative diagram.

\begin{center}
$%
\begin{array}
[c]{ccccc}
& X\times A & \overset{g^{\ast}\times A}{\longrightarrow} & Y^{A}\times A & \\
h_{X\times A} & \Downarrow & \searrow^{g} & \Downarrow & e_{Y}\\
& X\times A & \overset{g}{\rightarrow} & Y &
\end{array}
\nolinebreak$

Het adjunctive square for exponential endo-adjunction parsed as a coreflection
\end{center}

One might well ask: \textquotedblleft Why the special treatment since the
heteromorphisms are supposed to given (up to isomorphism) by the hom-pair
adjunctive square diagram of the adjunction representation
theorem?\textquotedblright\ The answer is that this is exactly what has been
derived. For the exponential adjunction, the hom-pair adjunctive square
diagram is as follows:

\begin{center}
$%
\begin{array}
[c]{ccccc}
& (X,X\times A) & \overset{\left(  f,g^{\ast}\times A\right)  }%
{\longrightarrow} & \left(  Y^{A},Y^{A}\times A\right)  & \\
\left(  \eta_{X},1_{X\times A}\right)  & \downarrow & \searrow^{(f,g)} &
\downarrow & \left(  1_{Y^{A}},\varepsilon_{Y}\right) \\
& \left(  \left(  X\times A\right)  ^{A},X\times A\right)  & \overset{(\left(
f^{\ast}\right)  ^{A},g)}{\longrightarrow} & \left(  Y^{A},Y\right)  &
\end{array}
\nolinebreak$

Hom-pair adjunctive square for exponential adjunction
\end{center}

The left component of each pair is the `same' as the het adjunctive square for
the adjunction parsed as a reflection (see the previous diagram for the
reflection) modulo the point that the codomains of the vertical maps are taken
as $\mathbf{APower}$ so they become heteromorphisms (e.g., $\eta_{X}$ becomes
$h_{X}$). Dually, the right component of each pair is the `same' as the het
adjunctive square for the adjunction parsed as a coreflection (see the
previous diagram) modulo the point that the domains of the vertical maps are
taken as $\mathbf{AProd}$ so they become heteromorphisms (e.g., $\varepsilon
_{Y}$ becomes $e_{Y}$).

\section{Concluding Remarks}

This paper inevitably has two themes. The main theme is showing how adjoint
functors arise from the representations within two categories of the
heteromorphisms between the categories. But the logically prior theme is
showing that heteromorphisms can be rigorously treated as part of category
theory---rather than just as stray chimeras roaming in the wilds of
mathematical practice.

Taking the logically prior theme first, category theory has always been
presented as embodying the idea of grouping mathematical objects of a certain
sort together with their appropriate morphisms in a \textquotedblleft
category.\textquotedblright\ In some respects, this homomorphic theme became
the leading theme just as in Felix Klein's Erlanger Program where geometries
were characterized by the invariants of a specified class of transformations.
Indeed, in their founding paper, Eilenberg and MacLane noted that category
theory \textquotedblleft may be regarded as a continuation of the Klein
Erlanger Program, in the sense that a geometrical space with its group of
transformations is generalized to a category with its algebra of
mappings.\textquotedblright\ \cite[p. 237]{eilmac:gte} Hence the whole concept
of a \textquotedblleft heteromorphism\textquotedblright\ between objects of
different categories has seemed like a cross-species hybrid that is
out-of-place and running against the spirit of the enterprise. Functors were
defined to handle all the external relations between categories so
object-to-object inter-category morphisms had no `official' role.

At the outset of this paper, a number of testimonials were quoted about the
centrality of adjunctions in category theory. Once heteromorphisms were
rigorously treated using het-bifunctor (in analogy to treating homomorphisms
with hom-bifunctors), it quickly became clear that an adjunction between two
categories was closely related to the heteromorphisms between the objects of
the two categories. Our main theme is that left and right adjoints arise as
the left and right representations within the categories of the
heteromorphisms between the categories. Given the importance of adjoints, this
made an argument for taking heteromorphisms \textquotedblleft out of the
closet\textquotedblright\ and recognizing them as part of the conceptual
family of category theory.\footnote{Indeed, one might ask why it has taken so
long for heteromorphisms to be formally recognized? The reticence seems to
come from heteromorphisms not fitting into the Erlanger-style theme of
emphasizing homomorphisms between mathematical objects of the same category.}\ 

When heteromorphisms (and the related concepts such as het natural
transformations) are accepted as a part of category theory---as is argued
here, then it does imply some \textquotedblleft broadening\textquotedblright%
\ of that Erlanger-style theme. Mathematical objects in different categories
may nevertheless have \textit{partial} similarities of structure which would
be expressed by heteromorphisms between the categories. Adjoints arise from
universals that represent and internalize those structural external
relationships within each of the categories. And that may help explain why
adjoints have emerged as the principal lens to focus on what is important in
mathematics. Structures are important \cite{ell:ae} that have within them
universal models representing the relationships with external entities of a
different kind.

\end{document}